\documentclass[11pt]{article}
\usepackage{amsfonts}
\usepackage{amssymb}
\usepackage{amsmath}
\usepackage{color}
=0pt
\def\sqr#1#2{{\vcenter{\vbox{\hrule height.#2pt
              \hbox{\vrule width.#2pt height#1pt \kern#1pt \vrule
width.#2pt}
              \hrule height.#2pt}}}}
\def\signed #1{{\unskip\nobreak\hfil\penalty50
              \hskip2em\hbox{}\nobreak\hfil#1
              \parfillskip=0pt \finalhyphendemerits=0 \par}}
\def\endpf{\signed {$\sqr69$}}
\def\dbE{{\mathbb{E}}}
\def\dbF{{\mathbb{F}}}

\def\dbN{{\mathbb{N}}}

\def\dbP{{\mathbb{P}}}

\def\dbR{{\mathbb{R}}}

\def\Om{\Omega}
\def\om{\omega}
\def\b{\beta}

\def\l{\lambda}

\def\si{\sigma}

\def\f{\varphi}

\def\3n{\negthinspace \negthinspace \negthinspace }
\def\2n{\negthinspace \negthinspace }
\def\1n{\negthinspace }
\def\ns{\noalign{\smallskip} }

\def\ds{\displaystyle}
\def\G{\Gamma}
\def\D{\Delta}

\def\L{\Lambda}

\def\cC{{\cal C}}

\def\cF{{\cal F}}

\def\cJ{{\cal J}}

\def\cL{{\cal L}}
\def\cM{{\cal M}}

\def\cX{{\cal X}}

\def\mE{{\mathbb{E}}}

\def\no{\noindent}

\def\bs{\bigskip}
\def\q{\quad}
\def\qq{\qquad}
\def\hb{\hbox}

\def\esssup{\mathop{\rm esssup}}

\def\pa{\partial}

\def\wt{\widetilde}
\def\cd{\cdot}
\def\cds{\cdots}

\def\as{\hbox{\rm a.s.{ }}}

\def\({\Big (}
\def\){\Big )}
\def\[{\Big[}
\def\]{\Big]}

\def\={\buildrel \triangle \over =}

\def\resp{{\it resp. }}

\def\be{\begin{equation}}
\def\bel{\begin{equation}\label}
\def\ee{\end{equation}}
\def\bea{\begin{eqnarray}}
\def\eea{\end{eqnarray}}
\def\bt{\begin{theorem}}
\def\et{\end{theorem}}
\def\bc{\begin{corollary}}
\def\ec{\end{corollary}}
\def\bl{\begin{lemma}}
\def\el{\end{lemma}}
\def\bp{\begin{proposition}}
\def\ep{\end{proposition}}
\def\br{\begin{remark}}
\def\er{\end{remark}}
\def\ba{\begin{array}}
\def\ea{\end{array}}
\def\bd{\begin{definition}}
\def\ed{\end{definition}}

\newtheorem{lemma}{Lemma}[section]
\newtheorem{remark}{Remark}[section]

\newtheorem{theorem}{Theorem}[section]
\newtheorem{corollary}{Corollary}[section]

\newtheorem{definition}{Definition}[section]
\newtheorem{proposition}{Proposition}[section]

\makeatletter
   
   \@addtoreset{equation}{section}
\makeatother

\begin{document}

\title{\bf Stochastic Well-posed Systems and Well-posedness
of Some Stochastic Partial Differential
Equations with Boundary Control and Observation}

\author{Qi L\"u\thanks{School of Mathematics,
Sichuan University, Chengdu, 610064, China ({\tt
luqi59@163.com}). This work is supported by the
NSF of China under grants 11101070 and 11471231,
and the Advanced Grants NUMERIWAVES/FP7-246775
of the European Research Council Executive
Agency.} }

\date{}

\maketitle

\begin{abstract}
We generalize the concept ``well-posed linear
system" to stochastic linear control systems and
study some basic properties of such kind
systems. Under our generalized definition, we
show the well-posedness of the stochastic heat
equation and the stochastic Schr\"odinger
equation with suitable boundary control and
observation operators, respectively.
\end{abstract}

\bs

\no{\bf 2010 Mathematics Subject
Classification}.  Primary 93E03; Secondary
60H15.

\bs

\no{\bf Key Words}. stochastic well-posed linear
system, stochastic heat equation, stochastic
Schr\"odinger equation.

\section{Introduction}\label{s1}

Let $H$, $U$ and $\wt U$ be three Hilbert spaces
which are identified with their dual spaces. Let
$A$ be the generator of a $C_0$-semigroup
$\{S(t)\}_{t\geq 0}$ on $H$. Denote by $H_{-1}$
the completion of $H$ with respect to the norm
$|x|_{H_{-1}}\=|(\b I-A)^{-1}x|_H$, where
$\b\in\rho(A)$ is fixed. Let   $B \in
\cL(U,H_{-1})$ and $\cC\in \cL(D(A),\wt U)$.

Consider the following control system:
\begin{equation}\label{system1}
\left\{
\begin{array}{ll}\ds
\frac{dY(t)}{dt}= A Y(t)+
Bu(t) & \mbox{ in } (0,+\infty),\\
\ns\ds Y(0)=Y_0,\\
\ns\ds Z(t)=\cC Y(t) & \mbox{ in } (0,+\infty).
\end{array}
\right.
\end{equation}
Here $Y_0\in H$ and $u\in L^2(0,+\infty;U)$.

In \eqref{system1}, the expected state space is
$H$.  If $B\in \cL(U,H)$, one can easily show
that $Y(t)\in H$ for all $t\geq 0$ by the
classical theory of evolution equations(see, for
instance, \cite[Chapter 3]{CZ}).

If $\cC\in
\cL(H,\wt U)$, then the observation $\cC Y(t)$
makes sense for all $H$-valued function
$Y(\cd)$.


However, in practical control system,
it is very common that the control and
observation operators are unbounded with respect
to the state space. Typical examples are systems
governed by partial differential equations(PDEs) in which actuators and sensors act on
lower-dimensional hypersurfaces or on the
boundary of a spatial domain.

The unboundedness of the control and observation
operators leads to  substantial technical
difficulties even for the formulation of the
state space. For instance, when $B\in
\cL(U,H_{-1})$, it seems that the natural state
space for the dynamics of \eqref{system1} should
be $H_{-1}$. On the other hand, to handle the
unbounded observation operator $\cC$, the
natural state space  should be $D(A)$.

In \emph{the} deterministic setting, to overcome
the above gap between the state spaces, people
introduced the notions of the admissible control
operator and the admissible observation
operator(see \cite{Salamon1,Salamon2} for
example). Furthermore, to study the feedback
control problem, people introduced the notion of
``well-posed linear system", which satisfies,
roughly speaking, that the map from the input
space $L^2(0,T;U)$ to the output space
$L^2(0,T;\wt U)$ is
bounded(\cite{Salamon2,TW,Weiss1,Weiss2}).

The well-posed linear systems form a very
general class whose basic properties are rich
enough to develop a parallel theory for the
theory of control systems with bounded control
and observation operators, such as feedback
control, dynamic stabilization and
tracking/disturbance rejection. Furthermore, the
well-posed linear systems are quite general and
covers many control systems described by partial
differential equations with actuators and
sensors acting on lower-dimensional
hyper-surfaces or on the boundary of a spatial
domain, which is the starting point of the study
of such kind of systems.

As far as we know, they are amount of references
concerning the  well-posed deterministic linear
systems, both for abstract control system and
controlled PDEs  in the last three decades (see
\cite{Curtain,GL,GS,GZ1,GZ2,TW,Staffans,WST,WCG}
and the rich references therein). However,  to
our best knowledge, there is a blank field for
the well-posedness of the stochastic linear
systems. In this  paper, we generalize the
notion of well-posed linear system to the
stochastic context by providing a formulation of
the stochastic well-posed linear system and some
basic properties.

\vspace{0.15cm}

We first introduce some notations. Let
$(\Om,\cF,\dbF,\dbP)$ be a complete filtered
probability space with a filtration
$\dbF=\{\cF_t\}_{t\ge0}$, on which a
one-dimensional standard Brownian motion
$\{W(t)\}_{t\ge0}$ is defined. Let $\cX$ be a
Banach space. For any $p\in [1,\infty)$, $t\geq
0$ and a sub-$\si$-algebra $\cM$ of $\cF$,
denote by $L_{\cM}^p(\Omega;\cX)$(\resp
$L_{\cM}^\infty(\Omega;\cX)$) the set of all
$\cM$-measurable ($\cX$-valued) random variables
$\xi:\Omega\to \cX$ with
$\mathbb{E}|\xi|_{\cX}^p < \infty$(\resp
$\esssup_{\om\in\Om}|\xi|_{\cX}<\infty$). Next,
for any $p,q\in[1,\infty)$ and $T\in
(0,+\infty]$, put
$$
\begin{array}{ll}\ds
L^p_\dbF(0,T;L^q(\Omega;\cX))\= \Big\{ \f :(0,T)
\times \Omega\!\to \cX \,\Bigm|\, \f(\cd)\hb{ is
$\dbF$-adapted   and
}\int_0^T\(\dbE|\f(t)|_\cX^q\)^{\frac{p}
{q}}dt<\infty\Big\},\\
\ns\ds L^\infty_\dbF(0,T;L^q(\Omega;\cX))\=
\Big\{ \f :(0,T) \times \Omega\!\to \cX
\,\Bigm|\, \f(\cd)\hb{ is $\dbF$-adapted   and
}\esssup_{t\in
[0,T]}\dbE|\f(t)|_\cX^q<\infty\Big\},\\
\ns\ds L^p_\dbF(0,T;L^\infty(\Omega;\cX))\=
\Big\{ \f :(0,T)\! \times \!\Omega\!\to \!\cX
\,\Bigm|\, \f(\cd)\hb{ is $\dbF$-adapted   and
}\int_0^T\!\!\big[\esssup_{\om\in \Om}
|\f(t)|_\cX\big]^pdt \!<\!\infty\Big\},\\
\ns\ds L^\infty_\dbF(0,T;\cX)\= \Big\{\f
:(0,T)\! \times \!\Omega\!\to \!\cX \,\Bigm|\,
\f(\cd)\hb{ is $\dbF$-adapted   and
}\esssup_{(t,\om)\in [0,T]\times\Om} |\f(t)|_\cX
<\infty\Big\},\\
\ns\ds L^p_{Br}(0,T;\cX)\= \Big\{ \f :(0,T) \to
\cX \,\Bigm|\, \f(\cd)\hb{ is Borel measurable
and }\int_0^T\(|\f(t)|_\cX^p\)dt<\infty\Big\}.
\end{array}
$$
When $p=q$, we simply denote
$L^p_\dbF(0,T;L^p(\Omega;\cX))$ by
$L^p_\dbF(0,T;\cX)$. Further, for any $p\in
[1,\infty)$, set
$$
\begin{array}{ll}\ds
C_{\dbF}([0,T];L^{p}(\Omega;\cX))\=\Big\{\f:[0,T]\times\Omega\to
\cX \bigm| \f(\cd)\hb{ is $\dbF$-adapted  and
}\dbE\big(|\f(t)|_\cX^p\big)^{\frac{1}{p}}\mbox{
is continuous }\Big\}.
\end{array}
$$
All the above spaces are endowed with the
canonical norms.

Consider the following  stochastic control
system:
\begin{equation}\label{system2}
\left\{
\begin{array}{ll}\ds
dY(t)=\big[A+F_1(t)\big]Y(t)dt +
Bu(t) dt + F_2(t)Y(t) dW(t)& \mbox{ in } (0,+\infty),\\
\ns\ds Y(0)=Y_0,\\
\ns\ds Z(t)=\cC Y(t) & \mbox{ in } (0,+\infty).
\end{array}
\right.
\end{equation}
Here $Y_0\in L^2_{\cF_0}(\Om;H)$, $u\in
L^2_\dbF(0,T;U)$ and  $F_1,F_2\in
L^\infty_\dbF(0,+\infty;\cL(H))$.

\vspace{0.1cm}

We first give the definition of the mild
solution to \eqref{system2}.

\begin{definition}\label{def-mild}
An $H$-valued stochastic process $Y(\cd)$ is
called a mild solution to \eqref{system2} if
\begin{enumerate}
  \item $Y(\cd)\in C_\dbF([0,+\infty);L^2(\Omega;H))$;
  \item For any $t\in [0,+\infty)$,
$$
\begin{array}{ll}\ds
Y(t) \3n&\ds= S(t)Y_0 + \int_0^t S(t -
s)F_1(s)Y(s)ds + \int_0^t S(t - s)Bu(s)ds\\
\ns&\ds\q + \int_0^t S(t - s) F_2(s)Y(s)dW(s).
\end{array}
$$
\end{enumerate}
\end{definition}

In general, the stochastic convolution $\int_0^t
S(t - s) F_2(s)Y(s)dW(s)$ is no longer a
martingale. Then, we cannot apply  It\^{o}'s
formula to mild solutions directly on most
occasions. This will limit the way to establish
suitable energy estimate for stochastic partial
differential equations. To avoid this
restriction, the notion of weak solution is
introduced.

\begin{definition}\label{def-weak}
A process $Y(\cd)\in
C_\dbF([0,+\infty);L^2(\Omega;H))$ is called a
weak solution to \eqref{system2} if for all
$t\in [0,+\infty)$ and $\psi\in D(A^*)$,
$$
\begin{array}{ll}\ds
\big\langle Y(t), \psi \big\rangle_{H}\3n &\ds=
\big\langle Y_0, \psi \big\rangle_{H} + \int_0^t
\big\langle  Y(s),A^*\psi\big\rangle_{H} +
\int_0^t \big\langle F_1(s)Y(s), \psi
\big\rangle_{H} ds  + \int_0^t \big\langle Y(s),
B^*\psi \big\rangle_{H}
ds\\
\ns&\ds \q  + \int_0^t \big\langle F_2(s)Y(s),
\psi \big\rangle_{H}dW(s),\;\dbP\mbox{-}\as
\end{array}
$$
\end{definition}
\begin{remark}
If $B\in \cL(U,H)$, then it is well-known that a
mild solution to \eqref{system2} is also a weak
solution to \eqref{system2}(see \cite{Prato} for
example). However, when $B\in \cL(U,H_{-1})$, as
far as we know, there is no result for the
relationship between these two kinds of
solutions.
\end{remark}
The rest of the paper is organized as follows:
in Section \ref{sec-2}, we give the formulation
of stochastic well-posed linear system and some
basic properties about it. Sections \ref{sec-3}
and \ref{sec-4} are devoted to the study the
well-posedness of controlled stochastic heat
equations and Schr\"odinger equations,
respectively. In Section \ref{sec-5}, we give
some further comments and present some open
problems.

Please note that in order to present the key
idea in a simple way, we do not  pursue the full
technical generality in this paper.


\section{Formulation and basic properties of stochastic well-posed linear systems}
\label{sec-2}

In this section, we give the formulation of a
stochastic well-posed linear system and study
some of its basic properties. First, we recall
the notion of the admissible control operator,
which is first introduced in the context of
deterministic control systems.  Then we show the
existence and uniqueness of the mild and weak
solutions to \eqref{system1} when $B$ is an
admissible control operator. Next, we recall the
concept of the admissible observation operator.
At last, we present the definition of the
stochastic well-posed linear system.


\subsection{Admissible control
operator}\label{sec-2.1}


The concept of admissible control operator is
motivated by the study of the solution to the
deterministic control system \eqref{system1},
where, people would like to study the operator
$B$ for which all the mild solutions $Y$ to
\eqref{system1} belong to $C([0,+\infty);H)$.
Such  operator is called admissible. In this
paper, we will show that it is also a suitable
notion when studying the solution to stochastic
control system \eqref{system2}.

Let  $t\in [0,+\infty)$. We define an operator
$\Phi_t \in
\cL(L^2_\dbF(0,+\infty;U),L^2_{\cF_t}(\Omega;H_{-1}))$
by
\begin{equation}\label{ad-con}
\Phi_t u = \int_0^t
S(t-s)B\big[\chi_{[0,t]}(s)u(s)\big]ds.
\end{equation}
\begin{definition}\label{def-ad-con}
The operator $B\in \cL(U,H_{-1})$ is called an
admissible control operator (for
$\{S(t)\}_{t\geq 0}$) if there is a $t_0>0$ such
that ${\rm Range}(\Phi_{t_0})\subset
L^2_{\cF_{t_0}}(\Omega;H)$.
\end{definition}
\begin{remark}
As far as we know,  the concept of admissible
control operator was first introduced in
\cite{HR} in the deterministic framework. Soon
afterwards, it was presented systematically as
an ingredient of the well-posed linear system in
\cite{Salamon2}.  Our definition of admissible
control operator for stochastic control systems
enjoys the same spirit(see Proposition
\ref{prop1} for the detail). However, we
formulate it as Definition \ref{def-ad-con} for
the convenience of the study of stochastic
control problems.
\end{remark}

\begin{remark}
Clearly, if $B$ is admissible, then in
\eqref{ad-con}, the integrand takes values in
$L^2_{\cF}(\Om;H_{-1})$, but the resultant
integral lies in $L^2_{\cF}(\Om;H)$, a dense
subspace of $L^2_{\cF}(\Om;H_{-1})$. As for the
deterministic case,  if $B\in\cL(U,H)$, then $B$
is admissible.
\end{remark}
\begin{remark}
People had proven the admissibility for some
control operators originated in controllability
problems of stochastic PDEs with boundary
control(see \cite{Lu2,Lu3} for example). In this
paper, we will prove that the control operators
in controlled  stochastic heat equations and
Schr\"odinger equations with suitable boundary
control are admissible.
\end{remark}
We have the following properties for the
admissible control operator.
\begin{proposition}\label{prop1}
The control operator $B\in \cL(U,H_{-1})$ is
admissible if and only if there is a constant
$C=C(t_0)>0$ such that for any $u\in
L^2_\dbF(0,+\infty;U)$,
\begin{equation}\label{def-ad-con1}
|\Phi_{t_0}u|_H \leq C\int_0^{t_0}
|u(s)|_U^2ds,\qq \dbP\mbox{-a.s.}
\end{equation}
\end{proposition}

{\it Proof}\,: The ``if" part is obvious. Let us
prove the ``only if" part. We do it by
contradiction argument. If \eqref{def-ad-con1}
were untrue, then there is a sequence
$\{u_n\}_{n=1}^\infty\subset
L^2_{Br}(0,+\infty;U)\subset
L^2_\dbF(0,+\infty;U)$ with
$|u_n|_{L^2_{Br}(0,+\infty;U)}=1$ such that
$$
\big|\Phi_{t_0} u_n\big|_{H}\geq n,
$$
that is,
\begin{equation}\label{12.25-eq5}
\big|\Phi_{t_0}
u_n\big|_{L^2_{\cF_{t_0}(\Om;H)}}\geq n.
\end{equation}

Let  $\l\in\rho(A)$ and $B_0=(\l I-A)^{-1}B$.
Since $B\in \cL(U,H_{-1})$, we have that
$B_0\in\cL(U,H)$ and
$$
\Phi_{t_0} u = (\l I-A)\int_0^{t_0}
S(t_0-s)B_0\big[\chi_{[0,t_0]}(s)u(s)\big]ds.
$$
This implies that $\Phi_{t_0}$ is closed.
According to the closed-graph theorem, we know
that $\Phi_{t_0}$ is bounded, i.e., there is a
constant $C(t_0)$ such that for any $u_n\in
L^2_{Br}(0,+\infty;U)$,
$$
\big|\Phi_{t_0}
u_n\big|_{L^2_{\cF_{t_0}}(\Om;H)}\leq
C(t_0)|u_n|_{L^2_\dbF(0,T;U)},
$$
a contradiction.
\endpf

\begin{proposition}\label{prop-ad0}
If  ${\rm Range}(\Phi_{t_0})\subset
L^2_{\cF_{t_0}}(\Omega;H)$ for a specific
$t_0\geq 0$, then for every $t> 0$, $\Phi_{t}\in
\cL(L^2_\dbF(0,+\infty;U),L^2_{\cF_t}(\Omega;H))$.
\end{proposition}

{\it Proof}\,: Let $u\in L^2_\dbF(0,+\infty;U)$.
It is easy to see that
$$
\begin{array}{ll}\ds
\Phi_{2t_0}u \3n&\ds
=\int_0^{2t_0}S(2t_0-s)B\big[\chi_{[0,2t_0]}(s)u(s)\big]ds \\
\ns&\ds =
S(t_0)\int_0^{t_0}S(t_0-s)B\big[\chi_{[0,t_0]}(s)u(s)\big]ds
+
\int_{t_0}^{2t_0}S(2t_0-s)B\big[\chi_{[t_0,2t_0]}(s)u(s)\big]ds
\\
\ns&\ds = S(t_0)\Phi_{t_0}u + \Phi_{t_0}\tilde
u.
\end{array}
$$
where $\tilde u(s) = u(t_0+s)$. According to
Proposition \ref{prop1}, we find that
$$
\begin{array}{ll}\ds
\mE\big|\Phi_{2t_0}u\big|_{H}^2 \3n&\ds=
\mE\big|S(t_0)\Phi_{t_0}u + \Phi_{t_0}\tilde
u\big|_H^2  \leq C\big(|u|_{L^2_\dbF(0,t_0;U)}^2
+ \mE|\tilde u|_{L^2(0,t_0;U)}^2\big)  \leq
C|u|_{L^2_\dbF(0,2t_0;U)}^2.
\end{array}
$$
This deduces that $\Phi_{2t_0}\subset
\cL(L^2_\dbF(0,+\infty;U),L^2_{\cF_{2t_0}}(\Omega;H))$.
By induction,  for all $n\in\dbN$,
$\Phi_{2^nt_0}\subset
\cL(L^2_\dbF(0,+\infty;U),L^2_{\cF_{2^nt_0}}(\Omega;H))$.

For any $t> 0$, there is a $n\in\dbN$ such that
$t\in (0,2^nt_0]$. For $u\in
L^2_\dbF(0,+\infty;U)$, let
$$
\tilde u(s)= \left\{
\begin{array}{ll}\ds
0, &\mbox{ if }s\in [0,2^nt_0-t),\\
\ns\ds u(s-t), &\mbox{ if } s\in
[2^nt_0-t,+\infty).
\end{array}
\right.
$$
Then, we have that
$$
\begin{array}{ll}\ds
\Phi_{t} u \3n&\ds=  \int_0^{t}
S(t-s)B\big[\chi_{[0,t]}(s)u(s)\big]ds\\
\ns&\ds =
\int_{2^nt_0-t}^{2^nt_0}S(2^nt_0-s)B\big[\chi_{[2^nt_0-t,2^nt_0]}(s)\tilde
u(s)\big]ds \\
\ns&\ds= \Phi_{2^nt_0}\tilde u.
\end{array}
$$
Hence, we get that
$$
\mE\big|\Phi_{t} u\big|_H^2
=\mE\big|\Phi_{2^nt_0}\tilde u\big|_H^2 \leq
C|\tilde u|_{L^2_\dbF(0,2^nt_0;U)}^2 =
C|u|_{L^2_\dbF(0,t;U)}^2,
$$
which concludes that $\Phi_t\in
\cL(L^2_\dbF(0,+\infty;U),L^2_{\cF_t}(\Omega;H))$.
\endpf

\begin{proposition}\label{prop-ad1}
If $B$ is an admissible control operator, then
the mapping
$$
\left\{
\begin{array}{ll} \Lambda:(0,+\infty)\times
L^2_\dbF(0,+\infty;U)\to L^2_{\cF}(\Omega;H),
\\
\ns\ds \L(t,u)=\Phi_t u,
\end{array}
\right.
$$
is continuous.
\end{proposition}

{\it Proof}\,: From Proposition \ref{prop-ad0},
we know that for any $t\in (0,+\infty)$,
$\L(t,\cd)$ is a continuous map from
$L^2_\dbF(0,+\infty;U)$ to
$L^2_{\cF_{t}}(\Omega;H)\subset
L^2_{\cF}(\Omega;H)$.

Next, we prove the continuity of $\L(t, u)$ with
respect to  $t$. Let $u\in
L^2_\dbF(0,+\infty;U)$ be fixed and put
$f(t)=\Phi_t u$. Let $0<t_1<t_2$. By Lebesgue's
dominated convergence theorem, we find that
$$
\begin{array}{ll}\ds
\q\lim_{t_2\to t_1^+} \mE|f(t_2)-f(t_1)|_H^2
\\
\ns\ds \leq 2\lim_{t_2\to
t_1^+}\mE\Big|\int_{t_1}^{t_2}S(t_2-s)\big[\chi_{[0,t_2]}(s)u(s)\big]ds
\Big|_H^2 \\
\ns\ds \q + 2\lim_{t_2\to
t_1^+}\mE\Big|\int_{0}^{t_1}\big[S(t_2-t_1)-I\big]S(t_1-s)\big[\chi_{[0,t_1]}(s)u(s)\big]ds
\Big|_H^2
\\
\ns\ds =0.
\end{array}
$$
This shows that $f(\cd)$ is right continuous.
Similarly, we can show that $f(\cd)$ is left
continuous. Hence, $f(\cd)$ is continuous.

The joint continuity of $\L(\cd,\cd)$ follows
easily from the fact that
$$
\Phi_t u - \Phi_s v = \Phi_t(u-v) +
(\Phi_t-\Phi_s)v.
$$
\endpf


\subsection{The existence and uniqueness of the mild and weak solution to \eqref{system2}}
\label{sec-2.2}

From Proposition \ref{prop-ad1}, we can prove
the following result.

\begin{theorem}\label{th-well-con1}
Let $B$ be an admissible control operator. Then
the equation \eqref{system2} admits a unique
mild solution $Y(\cd)\in
C_\dbF([0,+\infty);L^2(\Omega;H))$. Moreover,
for any $T>0$, there is a constant $C(T)>0$ such
that
\begin{equation}\label{th-well-con1-eq1}
|Y|_{C_\dbF([0,T];L^2(\Omega;H))} \leq
C(T)\big(|Y_0|_{L^2_{\cF_0}(\Om;H)}+
|u|_{L^2_\dbF(0,+\infty;U)} \big).
\end{equation}
\end{theorem}

Once we assume that $B$ is admissible, the proof
of Theorem \ref{th-well-con1} is very similar to
the case that $B\in \cL(U,H)$, which can be
found in \cite[Chapter 6]{Prato}. We give it
here not only for completeness but also for
presenting how we utilize the fact that $B$ is
admissible.

{\it Proof of Theorem \ref{th-well-con1}}\,: Let
$f,g\in L^2_\dbF(0,+\infty;H)$ . We claim that
\begin{equation}\label{12.25-eq3}
\begin{array}{ll}\ds
Y(\cd) \3n&\ds\= S(\cd)Y_0 + \int_0^\cd
S(t-s)f(s)ds +
\int_0^\cd S(t-s)Bu(s)ds\\
\ns&\ds\q + \int_0^\cd S(t-s)g(s)dW(s)\in
C_\dbF([0,+\infty);L^2(\Omega;H)).
\end{array}
\end{equation}

First,  we have that
\begin{equation}\label{12.25-eq1}
\begin{array}{ll}\ds
\q|Y(t)|_{L^2_{\cF_t}(\Omega;H)} \\
\ns\ds = \Big|S(t)Y_0 + \int_0^t S(t-s)f(s)ds +
\int_0^t S(t-s)Bu(s)ds + \int_0^t
S(t-s)g(s)dW(s)\Big|_{L^2_{\cF_t}(\Omega;H)}
\\
\ns \ds \leq
\big|S(\cd)Y_0\big|_{L^2_{\cF_0}(\Omega;H)}
 + \Big|\int_0^t
S(t - s)f(s)ds\Big|_{L^2_{\cF_t}(\Omega;H)} +
\Big|\int_0^t S(t-s)Bu(s)ds\Big|_{L^2_{\cF_t}(\Omega;H)}\\
\ns\ds \q  + \Big|\int_0^t S(t -
s)g(s)dW(s)\Big|_{L^2_{\cF_t}(\Omega;H)}
\\
\ns\ds \leq
C(t)\big(|Y_0|_{L^2_{\cF_0}(\Omega;H)} + |
u|_{L^2_\dbF(0,t;U)} + |f|_{L^2_\dbF(0,t;H)}
+|g|_{L^2_\dbF(0,t;H)}\big).
\end{array}
\end{equation}
Therefore, we find that
\begin{equation}\label{12.25-eq2}
\begin{array}{ll}\ds
|Y(\cd)|_{L^\infty_\dbF(0,T;L^2(\Omega;H))} \leq
C(T)\big(|Y_0|_{L^2_{\cF_0}(\Omega;H)} +
|u|_{L^2_\dbF(0,T;U)} +  |f|_{L^2_\dbF(0,T;H)} +
|g|_{L^2_\dbF(0,T;H)}\big).
\end{array}
\end{equation}
Further, for $0\leq t_1 \leq t_2<\infty$,
\begin{equation}\label{12.12-eq3}
\begin{array}{ll}\ds
\q|Y(t_2)-Y(t_1)|_{L^2_{\cF_{t_2}}(\Omega;H)} \\
\ns\ds = \Big|\big[S(t_2) \!-\! S(t_1)\big]Y_0 +
\int_{t_1}^{t_2}\!\!S(t_2\!-s)f(s)ds +
\int_{t_1}^{t_2}\!\!S(t_2\!-s)Bu(s)ds+
\int_{t_1}^{t_2}\!\!S(t_2\!-s)g(s)dW(s)
\\
\ns\ds \q  +
\int_{0}^{t_1}\big[S(t_2-t_1)-I\big]S(t_1-s)f(s)ds+
\int_{0}^{t_1}\big[S(t_2-t_1)-I\big]S(t_1-s)Bu(s)ds
\\
\ns\ds \q  +
\int_{0}^{t_1}\big[S(t_2-t_1)-I\big]S(t_1-s)g(s)dW(s)\Big|_{L^2_{\cF_{t_2}}(\Omega;H)}
\\
\ns\ds \leq \big|\big[S(t_2-t_1)-I\big]
S(t_1)Y_0\big|_{L^2_{\cF_{t_2}}(\Omega;H)} +
\Big|\int_{t_1}^{t_2}S(t_2-s)f(s)ds\Big|_{L^2_{\cF_{t_2}}(\Omega;H)}
\\
\ns\ds \q +
\Big|\int_{t_1}^{t_2}S(t_2-s)Bu(s)ds\Big|_{L^2_{\cF_{t_2}}(\Omega;H)}
+ \Big|\int_{t_1}^{t_2} S(t_2-s)g(s)dW(s)
\Big|_{L^2_{\cF_{t_2}} (\Omega;H)}
\\
\ns\ds \q + \Big| \int_{0}^{t_1}\!\! S(t_1\! -\!
s)\big[S(t_2\! -\! t_1)\! -
I\big]f(s)ds\Big|_{L^2_{\cF_{t_2}}\!(\Omega;H)}+
\Big|\int_{0}^{t_1}\!\!S(t_1\!-\!
s)\big[S(t_2\!-\! t_1)\!-\!
I\big]Bu(s)ds\Big|_{L^2_{\cF_{t_2}}\!
(\Omega;H)}\\
\ns\ds \q +
\Big|\int_{0}^{t_1}S(t_1-s)\big[S(t_2-t_1)-I\big]g(s)dW(s)\Big|_{L^2_{\cF_{t_2}}\!(\Omega;H)}.
\end{array}
\end{equation}
Since
$$
\big|\big[S(t_2-t_1)-I\big]
S(t_1)Y_0\big|_{H}\leq C |Y_0|_{H},
$$
by Lebesgue's dominated convergence theorem, we
get that
\begin{equation}\label{12.12-eq4}
\lim_{t_1\to t_2^-}\big|\big[S(t_2-t_1)-I\big]
S(t_1)Y_0\big|_{L^2_{\cF_{t_2}}(\Omega;H)}=0.
\end{equation}
Furthermore,
\begin{equation}\label{12.12-eq5}
\begin{array}{ll}\ds
\q\lim_{t_1\to
t_2^-}\Big|\int_{t_1}^{t_2}S(t_2-s)f(s)ds\Big|_{L^2_{\cF_{t_2}}(\Omega;H)}
\\
\ns\ds \leq \lim_{t_1\to
t_2^-}\int_{t_1}^{t_2}\big|S(t_2-s)f(s)\big|_{L^2_{\cF_{t_2}}(\Omega;H)}ds
\\
\ns\ds \leq C \lim_{t_1\to
t_2^-}\int_{t_1}^{t_2} |f(s)
|_{L^2_{\cF_{t_2}}(\Omega;H)}ds =0
\end{array}
\end{equation}
and
\begin{equation}\label{12.12-eq5.1}
\begin{array}{ll}\ds
\q\lim_{t_1\to
t_2^-}\Big|\int_{t_1}^{t_2}S(t_2-s)Bu(s)ds\Big|_{L^2_{\cF_{t_2}}(\Omega;H)}
\\
\ns\ds \leq C(t_2) \lim_{t_1\to
t_2^-}\int_{t_1}^{t_2} |u(s)
|_{L^2_{\cF_{t_2}}(\Omega;U)}ds =0.
\end{array}
\end{equation}
Clearly,
$$
\Big|\int_{t_1}^{t_2}S(t_2-s)g
(s)dW(s)\Big|_{L^2_{\cF_{t_2}}(\Omega;H)} \leq
|g|_{L^2_\dbF(t_1,t_2;H)}.
$$
Since
$$
\(\int_{t_1}^{t_2} |g(s)|^2_{H}ds\)(\om) \leq
\(\int_{0}^{t_2} |g(s)|^2_{H}ds\)(\om), \q
\dbP\mbox{-}a.s.,
$$
utilizing Lebesgue's dominated convergence
theorem, we find that
\begin{equation}\label{12.12-eq6}
\begin{array}{ll}\ds
\lim_{t_1\to
t_2^-}\Big|\int_{t_1}^{t_2}S(t_2-s)g(s)dW(s)\Big|_{L^2_{\cF_{t_2}}(\Omega;H)}
\3n&\ds\leq C\lim_{t_1\to t_2^-}
\mE\(\int_{t_1}^{t_2} |g(s)|^2_{H}ds\)\\
\ns&\ds=C \mE\(\lim_{t_1\to
t_2^-}\int_{t_1}^{t_2} |g(s)|^2_{H}ds\)=0.
\end{array}
\end{equation}
Further, since
$$
\big| \big[S(t_2-t_1)-I\big]f (s)\big|_{H} \leq
C|f (s)|_{H},
$$
using Lebesgue's dominated convergence theorem
again, we find that
\begin{equation}\label{12.12-eq7}
\begin{array}{ll}\ds
\q\lim_{t_1\to t_2^-}\Big|
\int_{0}^{t_1}S(t_1-s)\big[S(t_2-t_1)-I\big]f
(s)ds\Big|_{L^2_{\cF_{t_2}}(\Omega;H)}
\\
\ns\ds \leq C\lim_{t_1\to
t_2^-}\int_{0}^{t_2}\big|
\big[S(t_2-t_1)-I\big]f
(s)\big|_{L^2_{\cF_{t_2}}(\Omega;H)}ds =0.
\end{array}
\end{equation}
Similarly, thanks to
$$
\Big| \int_{0}^{t_1}S(t_1-s)B
u(s)ds\Big|_{L^2_{\cF_{t_2}}(\Omega;H)}\leq
C|u|_{L^2_\dbF(0,T;U)},
$$
it follows from Lebesgue's dominated convergence
theorem that
\begin{equation}\label{12.12-eq7.1}
\begin{array}{ll}\ds
\q\lim_{t_1\to t_2^-}\Big|
\big[S(t_2-t_1)-I\big]\int_{0}^{t_1}S(t_1-s)Bu
(s)ds\Big|_{L^2_{\cF_{t_2}}(\Omega;H)} =0.
\end{array}
\end{equation}
Next,
\begin{equation}\label{12.12-eq8}
\begin{array}{ll}\ds
\q\Big|\int_{0}^{t_1}S(t_1-s)\big[S(t_2-t_1)-I\big]g(s)dW(s)\Big|_{L^2_{\cF_{t_2}}(\Omega;H)}
\\
\ns\ds \leq
C\mE\(\Big|\int_{0}^{t_2}\big|S(t_1-s)\big[S(t_2-t_1)-I\big]g(s)\big|_{H}^2ds\)
\\
\ns\ds \leq C
\mE\(\Big|\int_{0}^{t_2}\big|\big[S(t_2-t_1)-I\big]g(s)\big|_{H}^2ds.
\end{array}
\end{equation}
It is clear that
$$
\int_{0}^{t_2}\big|\big[S(t_2-t_1)-I\big]g(s)\big|_{H}^2ds
\leq C\int_{0}^{t_2} |g(s) |_{H}^2ds.
$$
Then, Lebesgue's dominated convergence theorem
together with \eqref{12.12-eq8}, implies that
\begin{equation}\label{12.12-eq9}
\begin{array}{ll}\ds
\q\lim_{t_1\to
t_2^-}\Big|\int_{0}^{t_1}S(t_1-s)\big[S(t_2-t_1)-I\big]g(s)dW(s)\Big|_{L^2_{\cF_{t_2}}(\Omega;H)}
\\
\ns\ds \leq C\lim_{t_1\to t_2^-}
\mE\(\Big|\int_{0}^{t_2}\big|\big[S(t_2-t_1)-I\big]g(s)\big|_{H}^2ds\)=0.
\end{array}
\end{equation}
From \eqref{12.12-eq3}--\eqref{12.12-eq7} and
\eqref{12.12-eq9}, we obtain that
$$
\lim_{t_1\to
t_2^-}|Y(t_2)-Y(t_1)|_{L^2_{\cF_{t_2}}(\Omega;H)}=0.
$$
Similarly, we can get that
$$
\lim_{t_2\to
t_1^+}|Y(t_2)-Y(t_1)|_{L^2_{\cF_{t_2}}(\Omega;H)}=0.
$$
Thus, we prove that $Y(\cd)\in
C_\dbF([0,+\infty);L^2(\Om;H))$.

Fix any $T_1\in [0,+\infty)$. Let us define a
map
$$
\cJ:C_\dbF([0,T_1];L^2(\Omega;H))\to
C_\dbF([0,T_1];L^2(\Omega;H))
$$
as follows. For any $X(\cd)\in
C_\dbF([0,T_1];L^2(\Omega;H))$,
$$
\begin{array}{ll}\ds
Y(t)=\cJ(X)(t)\3n&\ds\=S(t)Y_0 + \int_0^t
S(t-s)F_1(s)X(s)ds + \int_0^t S(t-s)Bu(s)ds \\
\ns&\ds \q+ \int_0^t S(t-s)F_2(s)X(s)dW(s).
\end{array}
$$

$\cJ$ is well-defined following
\eqref{12.25-eq3}. We claim that if $T_1$ is
small enough, then $\cJ$ is contractive. Indeed,
let $Y_j=\cJ(X_j)$($j=1,2$) and $T_2>T_1$ be
fixed. Then, by \eqref{12.25-eq2}, we can find a
constant $C(T_2)>0$ such that
$$
\begin{array}{ll}\ds
\q|Y_1-Y_2|_{C_\dbF([0,T_1];L^2(\Omega;H))} \\
\ns\ds \leq C(T_2) \big(|F_1X_1 - F_1
X_2|_{L^2_\dbF(0,T_1;H)} + |F_2X_1 - F_2
X_2|_{L^2_\dbF(0,T_1;H)}\big)\\
\ns \ds \leq C(T_2)
\big(|F_1|_{L^\infty_\dbF(0,T_1;H)} +
|F_2|_{L^\infty_\dbF(0,T_1;H)}\big)|X_1-X_2|_{L^2_\dbF(0,T_1;H)}\\
\ns\ds \leq C(T_2)\sqrt{T_1}
\big(|F_1|_{L^\infty_\dbF(0,T_1;H)} +
|F_2|_{L^\infty_\dbF(0,T_1;H)}\big)|X_1-X_2|_{C_\dbF([0,T_1];L^2(\Om;H))}.
\end{array}
$$
Thus, we know that $\cJ$ is contractive if $T_1
<\frac{1}{C(T_2)^2}$. By means of the Banach
fixed point theorem, $\cJ$ has a unique fixed
point $Y(\cdot)\in
C_{\dbF}([0,T_1];L^{2}(\Omega;H))$. It is clear
that $Y(\cd)$ is a mild solution of the
following equation:
\begin{equation}\label{12.25-eq7.1}
\left\{
\begin{array}{lll}
\ds dY(t) = \big[ AX(t) + F_1(t)X(t)\big] dt + Bu(t)dt +  F_2(t)Y(t)dW(t) & \mbox{ in } (0,T_1],\\
\ns\ds Y(0)=Y_0.
\end{array}
\right.
\end{equation}
By \eqref{12.25-eq2} again,  we find that
$$
\begin{array}{ll}\ds
\q \mE|Y(t)|^2_{H}\\
\ns \ds\le C(T_1) \(\mE|Y_0|^2_{H} +
|u|_{L^2_\dbF(0,T;U)}^2 +
\int_0^t|F_1(s)Y(s))|^2_{H}ds +
\int_0^t|F_2(s)Y(s)|^2_{H}ds\)\\
\ns \ds\leq C(T_1) \[\mE|Y_0|^2_{H} +
|u|_{L^2_\dbF(0,T;U)}^2 +
\big(|F_1|_{L^\infty_\dbF(0,T;H)}
+|F_2|_{L^\infty_\dbF(0,T;H)}
\big)\int_0^t|Y(s))|^2_{H}ds\].
\end{array}
$$
This, together with the Gronwall's inequality,
implies that
\begin{equation}\label{12.25-eq9}
\big|Y(\cdot)\big|_{C_{\dbF}([0,T_1];L^{2}(\Omega;H))}
\leq C(T_1)\big[|Y_0|_{L^2_{\cF_0}(\Omega;H)} +
|u|_{L^2_\dbF(0,T;U)}^2\big].
\end{equation}
Repeating the above argument, we obtain the mild
solution to the equation \eqref{system2}. The
uniqueness of such solution  is obvious. The
desired estimate \eqref{th-well-con1-eq1}
follows from \eqref{12.25-eq9}. This completes
the proof of Theorem \ref{th-well-con1}.
\endpf

\begin{remark}
By Theorem \ref{th-well-con1}, any
$$
(Y_0,u(\cd))\in  L^2_{\cF_0}(\Om;H) \times
L^2_\dbF(0,+\infty;U)
$$
uniquely determines the state trajectory
$Y(\cd\,;Y_0,u)$ in an explicit way.
\end{remark}

As we said before, to apply the It\^o formula,
we should consider the weak solution to
\eqref{system2}. Fortunately, we have the
following result.

\vspace{0.2cm}

\begin{proposition}\label{th3}
The mild solution to \eqref{system2} is also a
weak solution to \eqref{system2} and vice versa.
\end{proposition}

\vspace{0.1cm}

{\it Proof}\,: We first prove that a weak
solution to \eqref{system2} is also mild
solution.

Assume that $Y(\cd)$ is a weak solution to
\eqref{system2}. For any $\psi\in D(A^*)$ and
$r\in [0,+\infty)$, we have that
\begin{equation}\label{12.24-eq1}
\begin{array}{ll}\ds
\big\langle Y(r), \psi \big\rangle_{H}\3n &\ds=
\big\langle Y_0, \psi \big\rangle_{H} + \int_0^r
\big\langle Y(s), A^*\psi \big\rangle_{H} ds+
\int_0^r \big\langle F_1(s)Y(s), \psi
\big\rangle_{H}
ds\\
\ns&\ds \q + \int_0^r \big\langle u(s), B^*\psi
\big\rangle_{H} ds + \int_0^r \big\langle
F_2(s)Y(s), \psi
\big\rangle_{H}dW(s),\;\dbP\mbox{-}\as
\end{array}
\end{equation}
We choose $\psi = S^*(t-r)A^*\phi$ for some
$\phi\in D((A^*)^2)$ and $t\in [r,+\infty)$.
From \eqref{12.24-eq1}, we get that
\begin{equation}\label{12.24-eq2}
\3n\begin{array}{ll}\ds
\q\big\langle Y(r), S^*(t-r)A^*\phi\big\rangle_H \\
\ns\ds = \big\langle Y_0,
S^*(t-r)A^*\phi\big\rangle_H + \int_0^r
\big\langle Y(s), A^*S^*(t-r)A^*\phi
\big\rangle_H ds\\
\ns \ds \q + \int_0^r \big\langle F_1(s)Y(s),S^*
(t - r)A^* \phi \big\rangle_H ds+ \int_0^r
\big\langle u(s),B^*S^*(t-r)A^* \phi
\big\rangle_{H} ds\\
\ns\ds \q + \int_0^r \big\langle F_2(s)Y(s),S^*
(t -r)A^*\phi \big\rangle_HdW(s).
\end{array}
\end{equation}
Integrating \eqref{12.24-eq2} with respect to
$r$ from $0$ to $t$, we obtain that
\begin{equation}\label{12.24-eq3}
\begin{array}{ll}\ds
\q\int_0^t\big\langle Y(r),
S^*(t-r)A^*\phi\big\rangle_Hdr \\
\ns\ds = \int_0^t\big\langle Y_0,
S^*(t-r)A^*\phi\big\rangle_H + \int_0^t\int_0^r
\big\langle Y(s), A^*S^*(t-r)A^*\phi
\big\rangle_H dsdr
\\
\ns \ds \q + \int_0^t\int_0^r \big\langle
F_1(s)Y(s),A^*S^*(t-r)A^*\phi \big\rangle_H dsdr
+ \int_0^t\int_0^r \big\langle
u(s),B^*S^*(t-r)A^* \phi \big\rangle_{H}
dsdr \\
\ns\ds \q + \int_0^t\int_0^r \big\langle
F_2(s)Y(s),S^*(t-r)A^*\phi \big\rangle_HdW(s)
dr.
\end{array}
\end{equation}
By a direct computation, we see that
\begin{equation}\label{12.24-eq4}
\begin{array}{ll}\ds
\int_0^t\big\langle Y_0,
S^*(t-r)A^*\phi\big\rangle_Hdr =
\int_0^t\big\langle AS(t-r)Y_0,
\phi\big\rangle_Hdr  = \big\langle S(t)Y_0,
\phi\big\rangle_H - \big\langle
Y_0,\phi\big\rangle_H.
\end{array}
\end{equation}
By Fubini's theorem, we find that
\begin{equation}\label{12.24-eq5}
\begin{array}{ll}\ds
\q\int_0^t\int_0^r \big\langle Y(s),
A^*S^*(t-r)A^*\phi \big\rangle_H dsdr
\\
\ns\ds =\int_0^t\int_s^t \big\langle Y(s),
A^*S^*(t-r)A^*\phi \big\rangle_H drds
\\
\ns\ds =\int_0^t \big\langle Y(s),
S^*(t-s)A^*\phi \big\rangle_H ds - \int_0^t
\big\langle Y(s), A^*\phi \big\rangle_H ds,
\end{array}
\end{equation}
\begin{equation}\label{12.24-eq6}
\begin{array}{ll}\ds
\q\int_0^t\int_0^r \big\langle F_1(s)Y(s),
S^*(t-r)A^*\phi \big\rangle_H dsdr
\\
\ns\ds =\int_0^t\int_s^t \big\langle F_1(s)Y(s),
S^*(t-r)A^*\phi \big\rangle_H drds
\\
\ns\ds =\int_0^t \big\langle F_1(s)Y(s),
S^*(t-s)\phi \big\rangle_H ds - \int_0^t
\big\langle F_1(s)Y(s), \phi \big\rangle_H ds,
\end{array}
\end{equation}
and
\begin{equation}\label{12.24-eq6.1}
\begin{array}{ll}\ds
\q\int_0^t\int_0^r \big\langle u(s),
B^*S^*(t-r)A^*\phi \big\rangle_H dsdr
\\
\ns\ds =\int_0^t\int_s^t \big\langle u(s),
B^*S^*(t-r)A^*\phi \big\rangle_H drds
\\
\ns\ds =\int_0^t \big\langle u(s),
B^*S^*(t-s)\phi \big\rangle_H ds - \int_0^t
\big\langle u(s), B^*\phi \big\rangle_H ds.
\end{array}
\end{equation}
 By the stochastic Fubini theorem, we obtain
that
\begin{equation}\label{12.24-eq7}
\begin{array}{ll}\ds
\q\int_0^t\int_0^r \big\langle F_2(s)Y(s),
S^*(t-r)A^*\phi \big\rangle_H dW(s)dr
\\
\ns\ds =\int_0^t\int_s^t \big\langle F_2(s)Y(s),
S^*(t-r)A^*\phi\big\rangle_H drdW(s)
\\
\ns\ds =\int_0^t \big\langle F_2(s)Y(s),
S^*(t-s)\phi \big\rangle_H dW(s) - \int_0^t
\big\langle F_2(s)Y(s), \phi \big\rangle_HdW(s)
.
\end{array}
\end{equation}
From \eqref{12.24-eq3}--\eqref{12.24-eq7}, we
end up with
\begin{equation}\label{12.24-eq8}
\begin{array}{ll}\ds
\Big\langle  Y(t) -  S(t)Y_0  -  \int_0^t
S(t-s)F_1(s)Y(s)ds - \int_0^t
S(t-s)Bu(s)ds \\
\ns\ds  - \int_0^t  S(t-s) F_2(s)Y(s)dW(s), \phi
\Big\rangle_H =0.
\end{array}
\end{equation}
Since $D((A^*)^2)$ is dense in $H$, we know that
\eqref{12.24-eq8} also holds for any $\phi\in
H$, which concludes that $Y$ is a mild solution
to \eqref{system2}.

\vspace{0.3cm}

Next, we  prove that a mild solution to
\eqref{system2} is also a weak solution.

Assume that $Y$ is a mild solution to
\eqref{system2}. Then for any $\psi\in D(A^*)$
and $r\in [0,+\infty)$, we have that
\begin{equation}\label{12.24-eq9}
\begin{array}{ll}\ds
\Big\langle Y(r) - S(r)Y_0 -\int_0^r
S(r-s)F_1(s)Y(s)ds  - \int_0^r S(t-s)Bu(s)ds
\\
\ns\ds   -\int_0^r S(r-s)F_2(s)Y(s)dW(s),
A^*\psi \Big\rangle_H=0.
\end{array}
\end{equation}
Integrating \eqref{12.24-eq9} from $0$ to $t$
with respect to $r$, we find that
\begin{equation}\label{12.24-eq10}
\begin{array}{ll}\ds
\q \int_0^t\big\langle Y(r),
A^*\psi\big\rangle_H
\\
\ns\ds =\int_0^t \big\langle S(r)Y_0, A^*\psi
\big\rangle_H dr + \int_0^t\int_0^r \big\langle
S(r-s)F_1(s)Y(s),A^*\psi \big\rangle_H ds dr \\
\ns\ds \q  + \int_0^t\int_0^r \big\langle
S(r-s)Bu(s),A^*\psi \big\rangle_H ds dr+
\int_0^t \int_0^r \big\langle S(r-s)F_2(s)Y(s),
A^*\psi\big\rangle_H dW(s)dr.
\end{array}
\end{equation}
First, we have that
\begin{equation}\label{12.24-eq11}
\int_0^t \big\langle S(r)Y_0, A^*\psi
\big\rangle_H dr = \int_0^t \big\langle
AS(r)Y_0, \psi \big\rangle_H dr = \big\langle
S(t)Y_0, \psi \big\rangle_H - \big\langle  Y_0,
\psi \big\rangle_H.
\end{equation}
By Fubini's theorem, we get that
\begin{equation}\label{12.24-eq12}
\begin{array}{ll}\ds
\q\int_0^t\int_0^r \big\langle
S(r-s)F_1(s)Y(s),A^*\psi \big\rangle_H ds dr \\
\ns\ds =\int_0^t\int_s^t \big\langle
AS(r-s)F_1(s)Y(s),\psi \big\rangle_H drds
\\
\ns\ds = \int_0^t \big\langle
S(t-s)F_1(s)Y(s),\psi \big\rangle_H  ds -
\int_0^t \big\langle F_1(s)Y(s),\psi
\big\rangle_H ds
\end{array}
\end{equation}
and
\begin{equation}\label{12.24-eq12.1}
\begin{array}{ll}\ds
\q\int_0^t\int_0^r \big\langle
S(r-s)Bu(s),A^*\psi \big\rangle_H ds dr \\
\ns\ds =\int_0^t\int_s^t \big\langle
AS(r-s)Bu(s),\psi \big\rangle_H drds
\\
\ns\ds = \int_0^t \big\langle S(t-s)Bu(s),\psi
\big\rangle_H  ds - \int_0^t \big\langle
u(s),B^*\psi \big\rangle_H ds.
\end{array}
\end{equation}
Thanks to the stochastic Fubini theorem, we
obtain that
\begin{equation}\label{12.24-eq13}
\begin{array}{ll}\ds
\q\int_0^t\int_0^r \big\langle
S(r-s)F_2(s)Y(s),A^*\psi \big\rangle_H dW(s) dr \\
\ns\ds =\int_0^t\int_s^t \big\langle
AS(r-s)F_2(s)Y(s),\psi \big\rangle_HdW(s)dr
\\
\ns\ds = \int_0^t \big\langle
S(t-s)F_2(s)Y(s),\psi \big\rangle_H dW(s) -
\int_0^t \big\langle F_2(s)Y(s),\psi
\big\rangle_HdW(s).
\end{array}
\end{equation}
From \eqref{12.24-eq10}--\eqref{12.24-eq13}, we
have that
\begin{equation}\label{12.24-eq14}
\begin{array}{ll}\ds
\big\langle Y(t),\psi\big\rangle_H \3n&\ds =
\big\langle Y_0,\psi\big\rangle_H + \int_0^t
\big\langle Y(s),A^*\psi \big\rangle_{H} ds +
\int_0^t \big\langle F_1(s)Y(s),\psi
\big\rangle_H ds \\
\ns&\ds\q +  \int_0^t \big\langle u(s),B^*\psi
\big\rangle_H ds+ \int_0^t \big\langle
F_2(s)Y(s),\psi \big\rangle_HdW(s),
\end{array}
\end{equation}
which shows that $Y$ is a weak solution to
\eqref{system2}.
\endpf


\subsection{Admissible observation
operator}\label{sec-2.3}


In this subsection, we study the admissible
observation operator. Let us first recall its
definition.

Define a family of operators $\{\Psi_t\}_{t\geq
0}$ from $D(A)$ to $L^2(0,+\infty;\wt U)$ as
follows:
\begin{equation}\label{12.25-eq4}
\big(\Psi_{t}\eta\big) (s)= \left\{
\begin{array}{ll}\ds
\cC S(s)\eta, &\mbox{ if } s\in [0,t],
\\
\ns\ds 0, &\mbox{ if } s\in (t,+\infty).
\end{array}
\right.
\end{equation}

\begin{definition}\label{def-ad-ob}
The operator $\cC \in\cL(D(A),\wt U)$ is called
an admissible observation operator for
$\{S(t)\}_{t\geq 0}$ if for some $t_0>0$,
$\Psi_{t_0}$ has a continuous extension to $H$.
\end{definition}

Although the definition of the admissible
observation operator comes from the
investigation of the deterministic control
system, we will show that it is also a suitable
notion in the study of stochastic control
system.

We first recall the following result.
\begin{proposition}\label{prop-adC2}
If $\cC \in\cL(D(A),\wt U)$ is admissible,  then
for every $t\geq 0$, $\Psi_t \in
\cL(H,L^2(0,+\infty;\wt U))$.
\end{proposition}
\begin{proposition}\label{prop3}
$\cC \in\cL(D(A),\wt U)$ is admissible if and
only if for every $t>0$, there is a constant
$C=C(t)>0$ such that
$$
\int_0^t |\cC S(s)\eta|_{\wt U}^2ds \leq
C(t)|\eta|_H^2,\q \forall\,\eta\in H.
$$
\end{proposition}

We refer the readers to \cite{Weiss1} for the
proofs of Propositions \ref{prop-adC2} and
\ref{prop3}.

In the stochastic context, we have to consider
the effect from the noise, i.e., the term
$\int_0^s S(s-r)Y(r)dW(r)$ and the fact that the
state space $L^2_{\cF_t}(\Om;H)$ depends on $t$.
Fortunately, we have the following result.

\begin{proposition}\label{prop2}
Let $\cC \in\cL(D(A),\wt U)$ be an admissible
observation operator for $\{S(t)\}_{t\geq 0}$.
Then for any $t>0$, there exists a constant
$C(t)
>0$ such that for any $Y_0\in
L^2_{\cF_0}(\Om;H)$, the solution to
\eqref{system2} with $u=0$ satisfies that
\begin{equation}\label{1.4-eq1.1}
\mE\int_{0}^{t}|\cC Y(s)|_{\wt U}^2ds \leq C(t)
\mE|Y_0|_H^2.
\end{equation}
\end{proposition}
\begin{remark}
If $Y(\cd)$ is a solution to a stochastic PDE
and $\cC$ is a boundary observation operator,
inequalities in the form of \eqref{1.4-eq1.1}
are usually called the hidden regularity of the
solution, i.e., it does not follow directly from
the classical trace theorem of Sobolev space. We
refer the readers to \cite{Lu4,Lu5,Lu3,Zhang}
for the hidden regularity for some stochastic
PDEs.
\end{remark}
{\it Proof of Proposition \ref{prop2}}\,: By the
closed-graph theorem, it is an easy matter to
see that there is a constant $C=C(t)>0$ such
that
\begin{equation}\label{12.25-eq6}
\int_0^t\big|\cC S(r)\eta\big|_H^2 dr \leq
C(t)|\eta|_H^2.
\end{equation}
Since
$$
Y(s) = S(s)Y_0 + \int_0^s S(s-r)F_1(r)Y(r)dr +
\int_0^s S(s-r)F_2(r)Y(r)dW(r),
$$
we have that
$$
\begin{array}{ll}\ds
\q\mE\int_{0}^{t}|\cC Y(s)|_{\wt U}^2ds \\
\ns\ds = \mE\int_{0}^{t}\Big|S(s)Y_0 + \int_0^s
S(s-r)F_1(r)Y(r)dr + \int_0^s
S(s-r)F_2(r)Y(r)dW(r)\Big|_{\wt U}^2 ds \\
\ns\ds \leq 3\mE\int_{0}^{t}\big|\cC
S(s)Y_0\big|_{\wt U}^2 ds +
3\mE\int_{0}^{t}\Big|\cC\int_0^s
S(s-r)F_1(r)Y(r)dr
\Big|_{\wt U}^2 ds \\
\ns\ds \q + 3\mE\int_{0}^{t}\Big|\cC\int_0^s
S(s-r)F_2(r)Y(r)dW(r)\Big|_{\wt U}^2 ds\\
\ns\ds \leq C(t)\(\mE\int_{0}^{t}\big|\cC
S(s)Y_0\big|_{\wt U}^2 ds +
\mE\int_{0}^{t}\big|\cC S(s-r)F_1(r)Y(r)
\big|_{\wt U}^2 ds + \mE\int_{0}^{t}\big|\cC
S(s-r)F_2(r)Y(r) \big|_{\wt U}^2 ds\).
\end{array}
$$
This, together with Proposition \ref{prop3},
implies that
$$
\begin{array}{ll}\ds
\q\mE\int_{0}^{t}|\cC Y(s)|_{\wt U}^2ds \\
\ns\ds   \leq C(t)\(\mE\big|Y_0\big|_{H}^2 ds +
\mE\int_{0}^{t}\big|F_1(s)Y(s) \big|_{H}^2 ds +
\mE\int_{0}^{t}\big|F_2(s)Y(s) \big|_{H}^2 ds\)\\
\ns\ds \leq C(t)\mE\big|Y_0\big|_{H}^2.
\end{array}
$$
\endpf

\begin{remark}
Obviously, every $\cC \in\cL(H,\wt U)$ is
admissible for $\{S(t)\}_{t\geq 0}$. If $\cC$ is
an admissible observation operator for
$\{S(t)\}_{t\geq 0}$, then we denote the
(unique) extension of $\Psi_{t}$ to
$L^2_{\cF_0}(\Omega;H)$ by the same symbol.
\end{remark}


\subsection{Stochastic well-posed linear
system}\label{sec-2.4}


We begin with the definition of the stochastic
well-posed linear system.

\begin{definition}\label{def-well}
Let $B\in \cL(U;H_{-1})$ be an admissible
control operator and $\cC\in \cL(D(A),\wt U)$ be
an admissible observation operator. We say
\eqref{system2} is well-posed if   there is a
$t_0>0$ and a $C(t_0)>0$ such that  for any
$Y_0\in L^2_{\cF_0}(\Om;H)$ and $u\in
L^2_\dbF(0,+\infty; U)$,
\begin{equation}\label{12.7-eq1}
\big|\cC Y\big|_{L^2_\dbF(0,t_0;\wt U)}\leq
C(t_0)\big(|Y_0|_{L^2_{\cF_0}(\Om;H)}+|u|_{L^2_\dbF(0,t_0;
U)}\big).
\end{equation}
\end{definition}
\begin{remark}
Although our definition of stochastic well-posed
linear system seems different from the classical
definition to the deterministic well-posed
linear system,  the spirit is the same, i.e.,
both the control and observation operators are
admissible and the map from the input to the
output is bounded.
\end{remark}
Similar to Propositions \ref{prop-ad0} and
\ref{prop-adC2}, we have the following result.
\begin{proposition}
If the system \eqref{system2} is well-posed,
then for any $t>0$, there is a constant $C(t)>0$
such that for any $Y_0\in L^2_{\cF_0}(\Om;H)$
and $u\in L^2_\dbF(0,+\infty; U)$,
\begin{equation}\label{12.7-eq2}
\big|\cC Y\big|_{L^2_\dbF(0,t;\wt U)}\leq
C(t)\big(|Y_0|_{L^2_{\cF_0}(\Om;H)}+|u|_{L^2_\dbF(0,t;
U)}\big).
\end{equation}
\end{proposition}
{\it Proof}\,: For any $u\in
L^2_\dbF(0,+\infty;U)$, by the uniqueness of the
solution to \eqref{system2}, we have that
$$
\chi_{[0,2t_0]}Y(\cd;0,u)=\chi_{[0,t_0]}Y\big(\cd;Y_0,\chi_{[0,t_0]}u\big)
+
\chi_{[t_0,2t_0]}Y\big(\cd;Y(t_0;Y_0,\chi_{[0,t_0]}u),\chi_{[t_0,2t_0]}u\big).
$$
Hence, we see that
\begin{equation}\label{12.25-eq8}
\begin{array}{ll}\ds
\q\big|\cC Y\big|_{L^2_\dbF(0,2t_0;\wt U)}\\
\ns\ds\leq \big|\cC
\big[\chi_{[0,t_0]}Y\big(\cd;Y_0,\chi_{[0,t_0]}u\big)\big]\big|_{L^2_\dbF(0,t_0;\wt
U)} + \big|\cC
\big[\chi_{[t_0,2t_0]}Y\big(\cd;Y(t_0;Y_0,\chi_{[0,t_0]}u),\chi_{[t_0,2t_0]}u\big)\big]\big|_{L^2_\dbF(t_0,2t_0;\wt
U)}.
\end{array}
\end{equation}
From \eqref{12.7-eq1}, we have that
\begin{equation}\label{12.25-eq10}
\big|\cC
\big[\chi_{[0,t_0]}Y\big(\cd;Y_0,\chi_{[0,t_0]}u\big)\big]\big|_{L^2_\dbF(0,t_0;\wt
U)} \leq
C(t_0)\big(|Y_0|_{L^2_{\cF_0}(\Om;H)}+|u|_{L^2_\dbF(0,t_0;
U)}\big).
\end{equation}
Next, put $\widehat
Y(t)=Y\big(t+t_0;Y(t_0;Y_0,\chi_{[0,t_0]}u),\chi_{[t_0,2t_0]}u\big)$.
Then we have that(see \cite[Chapter 4]{PR} for
the details)
\begin{equation}\label{12.25-eq11}
\begin{array}{ll}\ds
\q\big|\cC \big(\chi_{[t_0,2t_0]}\widehat
Y\big)\big|_{L^2_\dbF(t_0,2t_0;\wt U)}^2\\
\ns\ds = \mE\Big|\cC S(s)\widehat Y(t_0) +
\cC\int_{0}^{t_0}S(t_0-s)F_1(s)\widehat Y(s)ds +
\cC\int_{0}^{t_0}S(t_0-s)Bu(s+t_0)ds
\\
\ns\ds \qq+
\cC\int_{0}^{t_0}S(t_0-s)F_2(s)\widehat
Y(s)d\big[W(s+t_0)-W(t_0)\big]\Big|_{\wt U}^2 \\
\ns\ds \leq \!2 \mE\Big|\cC S(s)\widehat
Y(t_0)\! +\!
\cC\!\int_{0}^{t_0}\!\!S(t_0\!-\!s)F_1(s)\widehat
Y(s)ds\!+\!
\cC\!\int_{0}^{t_0}\!\!S(t_0\!-\!s)F_2(s)\widehat
Y(s)d\big[W(s\!+\!t_0)\!-\!W(t_0)\big]\Big|_{\wt U}^2 \\
\ns\ds \q + 2
\mE\Big|\cC\int_{0}^{t_0}S(t_0-s)Bu(s+t_0)ds\Big|_{\wt
U}^2.
\end{array}
\end{equation}
Since $\cC$ is an admissible observation
operator, we know that
\begin{equation}\label{12.25-eq12}
\begin{array}{ll}\ds
2 \mE\Big|\cC S(s)\widehat Y(t_0) \!+\!
\cC\!\int_{0}^{t_0}\!S(t_0\!-\!s)F_1(s)\widehat
Y(s)ds\!+\!
\cC\!\int_{0}^{t_0}\!\!S(t_0\!-\!s)F_2(s)\widehat
Y(s)d\big[W(s\!+\!t_0)\!-\!W(t_0)\big]\Big|_{\wt U}^2 \\
\ns\ds \leq C\(\mE|\widehat Y(t_0)|_H^2 +
\mE\int_{0}^{t_0}|\widehat Y(s)|_{H}^2 ds \).
\end{array}
\end{equation}
By Choosing $Y_0=0$ and noting that the system
\eqref{system2} is well-posed, from
\eqref{12.7-eq1}, we get that
\begin{equation}\label{12.25-eq13}
\mE\Big|\cC\int_{0}^{t_0}S(t_0-s)Bu(s)ds\Big|_{\wt
U}^2 \leq C(t_0)|u|_{L^2_\dbF(0,t_0;U)}^2.
\end{equation}
Hence, we have that for any $u\in
L^2_{Br}(0,+\infty;U)\subset
L^2_\dbF(0,+\infty;U)$,
\begin{equation}\label{12.25-eq14}
\begin{array}{ll}\ds
\Big|\cC\int_{0}^{t_0}S(t_0-s)Bu(s)ds\Big|_{\wt
U}^2 \3n&\ds =
\mE\Big|\cC\int_{0}^{t_0}S(t_0-s)Bu(s)ds\Big|_{\wt
U}^2\\
\ns&\ds \leq
C(t_0)\mE\int_0^{t_0}|u(s)|_{U}^2ds\\
\ns&\ds =C(t_0) \int_0^{t_0}|u(s)|_{U}^2ds.
\end{array}
\end{equation}
From \eqref{12.25-eq14}, we obtain that
\begin{equation}\label{12.25-eq15}
\mE\Big|\cC\int_{0}^{t_0}S(t_0-s)Bu(s+t_0)ds\Big|_{\wt
U}^2 \leq
C(t_0)\mE\int_0^{t_0}|u(s+t_0)|_{U}^2ds =
C(t_0)|u|_{L^2_\dbF(t_0,2t_0;U)}^2.
\end{equation}
According to
\eqref{12.25-eq8}--\eqref{12.25-eq12} and
\eqref{12.25-eq15}, we conclude that
\begin{equation}\label{12.25-eq16}
\big|\cC Y\big|_{L^2_\dbF(0,2t_0;\wt U)} \leq
C(2t_0)\big(|u|_{L^2_\dbF(0,2t_0;U)} +
\big|Y_0\big|_{L^2_{\cF_0}(\Om;H)} \big).
\end{equation}
By induction, we can prove that for any
$n\in\dbN$, there is a constant $C(n)>0$ such
that for any $Y_0\in L^2_{L^2_{\cF_0}(\Om;H)}$
and $u\in L^2_\dbF(0,+\infty;U)$,
\begin{equation}\label{12.17-eq3}
\big|\cC Y\big|_{L^2_\dbF(0,2^nt_0;\wt U)} \leq
C(n)\big(|u|_{L^2_\dbF(0,2^nt_0;U)} +
\big|Y_0\big|_{L^2_{\cF_0}(\Om;H)} \big).
\end{equation}

For any $t>0$, there is a $n\in\dbN$ such that
$0<t<2^n t_0$. For a given $u\in
L^2_\dbF(0,+\infty;U)$, put
$$
\tilde u(s) = \left\{
\begin{array}{ll}\ds
u(s), &\mbox{ if } s\in [0,t],\\
\ns\ds 0, &\mbox{ if } s\in (t,+\infty).
\end{array}
\right.
$$
According to \eqref{12.17-eq3}, we get that
\begin{equation}\label{12.17-eq4}
\begin{array}{ll}\ds
\big|\cC Y\big|_{L^2_\dbF(0,t;\wt U)}
\3n&\ds\leq C(n)\big(|\tilde
u|_{L^2_\dbF(0,2^nt_0;U)} +
\big|Y_0\big|_{L^2_{\cF_0}(\Om;H)} \big)\\
\ns&\ds\leq C(n)\big(|u|_{L^2_\dbF(0,t;U)} +
\big|Y_0\big|_{L^2_{\cF_0}(\Om;H)} \big).
\end{array}
\end{equation}
\endpf


\section{Well-posedness of controlled stochastic heat equation}
\label{sec-3}

In this section, we study a stochastic heat
equation with boundary control and observation.
We show that it is a stochastic well-posed
linear system.

Let $G\subset\dbR^d$($d\in\dbN$) be a bounded
domain with the $C^2$ boundary $\G$. Consider
the following stochastic heat equation:
\begin{equation}\label{H-system1}
\left\{
\begin{array}{ll}\ds
dy-\D y dt = aydt + bydW(t) &\mbox{ in }
G\times (0,+\infty),\\
\ns\ds \frac{\pa y}{\pa\nu}=u &\mbox{ on }\G \times (0,+\infty),\\
\ns\ds y(0)=y_0 &\mbox{ in } G,\\
\ns\ds z =y &\mbox{ on }\G \times (0,+\infty).
\end{array}
\right.
\end{equation}
Here $a,b\in
L^\infty_\dbF(0,+\infty;L^\infty(G))$ and
$y_0\in L^2_{\cF_0}(\Om;L^2(G))$.
\begin{theorem}\label{th2}
With the choice that $H=L^2(G)$,
$U=H^{-\frac{1}{2}}(\G)$ and $\wt
U=H^{\frac{1}{2}}(\G)$, the system
\eqref{H-system1} is well-posed, i.e., for any
$T>0$, there is a constant $C=C(T)>0$ such that
for any $y_0\in L^2_{\cF_0}(\Om;H)$ and $u\in
L^2_\dbF(0,T;U)$, there is a unique mild (also
weak) solution $y\in C_\dbF([0,T];L^2(\Om;H))$
to \eqref{H-system1} such that
\begin{equation}\label{th2-eq1}
|y|_{C_\dbF([0,T];L^2(\Om;H))} +
|z|_{L^2_\dbF(0,T;\wt U)} \leq C(T)\big(
|y_0|_{L^2_{\cF_0}(\Om;H)} +
|u|_{L^2_\dbF(0,T;U)} \big).
\end{equation}
\end{theorem}

{\it Proof}\,: We divide the proof into three
steps.

\vspace{0.12cm}

{\bf Step 1}. We first handle the case that
$u\in
C_\dbF^1([0,T];L^2(\Om;H^{-\frac{1}{2}}(\G)))$.

Consider the following elliptic equation:
$$
\left\{
\begin{array}{ll}\ds
\D v(t,\om) = 0 &\mbox{ in } G,\\
\ns\ds \frac{\pa v(t,\om)}{\pa\nu}=u(t,\om)
&\mbox{ on } \G_0.
\end{array}
\right.
$$
We claim that
\begin{equation}\label{12.18-eq3}
v\in C^1_\dbF([0,T];L^2(\Om;H^1(G))).
\end{equation}
Indeed, from the classical theory of elliptic
equations with Neumann boundary condition(see
\cite[Chapter 3]{LU} for example), we have that
\begin{equation}\label{12.18-eq4}
|v|^2_{H^{1}(G)} \leq C
|u|^2_{H^{-\frac{1}{2}}(\G)},\q\forall\,t\in
[0,T],\;\dbP\mbox{-a.s.}
\end{equation}
This, together with the fact that $u(t)\in
L^2_{\cF_t}(\Om;H^{-\frac{1}{2}}(\G))$, implies
that $v(t)\in L^2_{\cF_t}(\Om;H^{1}(G))$.
Furthermore, it follows from \eqref{12.18-eq4}
that
\begin{equation}\label{12.18-eq5}
\mE|v(s)-v(t)|^2_{H^{1}(G)} \leq C
\mE|u(s)-u(t)|^2_{H^{-\frac{1}{2}}(\G)}.
\end{equation}
Since $u\in
C_\dbF^1([0,T];L^2(\Om;H^{-\frac{1}{2}}(\G)))$,
we find from \eqref{12.18-eq5} that $v\in
C_\dbF^1([0,T];L^2(\Om;H^{1}(\G)))$.

Consider the following stochastic heat equation:
\begin{equation}\label{H-system2}
\left\{
\begin{array}{ll}\ds
d\tilde y-\D \tilde y dt = a\tilde ydt + (av
-v_t)dt + b\tilde ydW(t)+ bvdW(t) &\mbox{ in }
G\times (0,+\infty),\\
\ns\ds \frac{\pa \tilde y}{\pa\nu}=0 &\mbox{ on }\G \times (0,+\infty),\\
\ns\ds \tilde y(0)=y_0-v(0) &\mbox{ in } G.
\end{array}
\right.
\end{equation}
According to the classical well-posedness result
for stochastic heat equation(see \cite[Chapter
6]{Prato} for example), we know that
\eqref{H-system2} admits a unique mild solution
$\tilde y\in C_\dbF([0,T];L^2(\Om;L^2(G)))\cap
L^2_\dbF(0,T;H^1(G))$, which is also a weak
solution to \eqref{H-system2}. Let $y=\tilde y +
v$. Then, it is easy to see that $y$ is a
solution to \eqref{H-system1}.

{\bf Step 2}. In this step, we establish an
energy estimate for the solution to
\eqref{H-system1}.

By It\^o's formula, we have that
\begin{equation}\label{12.18-eq1}
\begin{array}{ll}\ds
\q\mE\int_G |y(t)|^2dx + \mE\int_0^t\int_G
|\nabla
y|^2dxds\\
\ns\ds = \mE\int_G |y(0)|^2dx + 2\mE\int_0^t
\langle u,y
\rangle_{H^{-\frac{1}{2}}(\G),H^{\frac{1}{2}}(\G)}
ds + 2\mE\int_0^t\int_G a y^2dxds +
\mE\int_0^t\int_G b^2 y^2dxds.
\end{array}
\end{equation}
Thanks to the classical trace theorem in Sobolev
space, we get that
$$
\begin{array}{ll}\ds
\q\mE\int_0^t \langle u,y
\rangle_{H^{-\frac{1}{2}}(\G),H^{\frac{1}{2}}(\G)}
ds\\
\ns\ds \leq \mE\int_0^t
|u|_{H^{-\frac{1}{2}}(\G)}|y|_{H^{\frac{1}{2}}(\G)}
ds \leq C\mE\int_0^t
|u|_{H^{-\frac{1}{2}}(\G)}|y|_{H^{1}(G)} ds\\
\ns \ds \leq 2C \mE\int_0^t
|u|^2_{H^{-\frac{1}{2}}(\G)}ds +
\frac{1}{4}\mE\int_0^t|y|^2_{H^{1}(G)} ds\\
\ns \ds \leq 2C\mE\int_0^t
|u|^2_{H^{-\frac{1}{2}}(\G)}ds +
\frac{1}{4}\mE\int_0^t\big(|\nabla y|^2 +
y^2\big) dxds.
\end{array}
$$
This, together with \eqref{12.18-eq1}, implies
that
\begin{equation}\label{12.18-eq2}
\begin{array}{ll}\ds
\q\mE\int_G |y(t)|^2dx +
\frac{1}{2}\mE\int_0^t\int_G |\nabla
y|^2dxds\\
\ns\ds \leq \mE\int_G |y(0)|^2dx + 4C
\mE\int_0^t |u|^2_{H^{-\frac{1}{2}}(\G)}ds +
\mE\int_0^t\int_G (2a + b^2 + 1) y^2dxds\\
\ns\ds \leq \mE\int_G |y(0)|^2dx \!+\! 4C
\mE\int_0^t |u|^2_{H^{-\frac{1}{2}}(\G)}ds +\!
\big(2|a|_{L^\infty_\dbF(0,T;L^\infty(G))} \!+\!
|b|^2_{L^\infty_\dbF(0,T;L^\infty(G))}\!+1\big)\mE\int_0^t\int_G
y^2dxds.
\end{array}
\end{equation}
It follows from Gronwall's inequality and
\eqref{12.18-eq2} that
\begin{equation}\label{12.18-eq2.1}
\mE\int_G |y(t)|^2dx +  \mE\int_0^T\int_G
|\nabla y|^2dxds \leq C\(\mE\int_G |y(0)|^2dx +
\mE\int_0^T |u|^2_{H^{-\frac{1}{2}}(\G)}ds\).
\end{equation}

{\bf Step 3}. In this step,  let us deal with
the case that $u\in
L^2_\dbF(0,T;H^{-\frac{1}{2}}(\G))$.

We can find a sequence $\{u_n\}_{n=1}^\infty
\subset
C^1_\dbF([0,T];L^2(\Om;H^{-\frac{1}{2}}(\G)))$
such that
\begin{equation}\label{12.18-eq6}
\lim_{n\to\infty} u_n = u \q\mbox{ in
}L^2_\dbF(0,T;H^{-\frac{1}{2}}(\G)).
\end{equation}
Denote by $y_n$ the solution to
\eqref{H-system1} with the initial datum $y_0$
and the control $u_n$. From \eqref{12.18-eq2.1}
and \eqref{12.18-eq6}, we know that
$\{y_n\}_{n=1}^\infty$ is a Cauchy sequence in
$C_\dbF([0,T];L^2(\Om;L^2(G)))\times
L^2_\dbF(0,T;H^1(G))$. Thus, there is a unique
$y\in C_\dbF([0,T];L^2(\Om;L^2(G)))\times
L^2_\dbF(0,T;H^1(G))$ such that
\begin{equation}\label{12.18-eq7}
\lim_{n\to\infty} y_n=y \q\mbox{ in
}\;C_\dbF([0,T];L^2(\Om;L^2(G)))\times
L^2_\dbF(0,T;H^1(G)).
\end{equation}
From the definition of $y_n$, we have that for
any $t\in [0,T]$ and $\eta\in H^1(G)$,
\begin{equation}\label{12.18-eq8}
\begin{array}{ll}\ds
\q\int_G y_n(t)\eta dx - \int_G y_n(0)\eta dx -
\int_0^t \langle u_n,\eta
\rangle_{H^{-\frac{1}{2}},H^{\frac{1}{2}}}ds +
\int_0^t\int_G \nabla y_n\nabla\eta dxds\\
\ns\ds = \int_0^t\int_G a y_n\nabla\eta dxds +
\int_0^t\int_G a y_n\nabla\eta
dxdW(s),\qq\dbP\mbox{-}\as
\end{array}
\end{equation}
Thanks to \eqref{12.18-eq7} and
\eqref{12.18-eq8}, we conclude that for any
$t\in [0,T]$,
\begin{equation}\label{12.18-eq9}
\begin{array}{ll}\ds
\q\int_G y(t)\eta dx - \int_G y(0)\eta dx -
\int_0^t \langle u,\eta
\rangle_{H^{-\frac{1}{2}},H^{\frac{1}{2}}}ds +
\int_0^t\int_G \nabla y\nabla\eta dxds\\
\ns\ds = \int_0^t\int_G a y\nabla\eta dxds +
\int_0^t\int_G a y\nabla\eta
dxdW(s),\qq\dbP\mbox{-}\as
\end{array}
\end{equation}
Hence, $y$ is a weak solution (also a mild
solution) to \eqref{H-system1} and satisfies
\eqref{th2-eq1}.
\endpf


\section{Well-posedness of controlled stochastic Schr\"odinger equation}
\label{sec-4}

This section is devoted to the study of a
stochastic Schr\"odinger equation with boundary
control and observation. We show that it is a
stochastic well-posed linear system.

In this section, we assume that $\dbF$ is the
natural filtration generated by the Brownian
motion $\{W(t)\}_{t\geq 0}$. Consider the
following stochastic Schr\"odinger equation:
\begin{equation}\label{S-system1}
\left\{
\begin{array}{ll}\ds
dy+i\D y dt = aydt + bydW(t) &\mbox{ in }
G\times (0,+\infty),\\
\ns\ds y=u &\mbox{ on }\G \times (0,+\infty),\\
\ns\ds y(0)=y_0 &\mbox{ in } G,\\
\ns\ds \f = -i\frac{\pa(-\D)^{-1}y}{\pa\nu}
&\mbox{ on }\G\times (0,+\infty).
\end{array}
\right.
\end{equation}
Here $y_0\in H^{-1}(G)$ and $a,b\in
L^\infty_\dbF(0,+\infty;W_0^{1,+\infty}(G))$.

Let $H=H^{-1}(G)$ and $U=\wt U=L^2(\G_0)$. We
have the following result:
\begin{theorem}\label{th1}
System \eqref{S-system1} is well-posed, i.e.,
for any $T>0$, there is a constant $C=C(T)>0$
such that for any $y_0\in L^2_{\cF_0}(\Om;H)$
and $u\in L^2_\dbF(0,T;U)$, there is a unique
solution $y\in C_\dbF([0,T];$ $L^2(\Om;H))$ to
\eqref{S-system1} such that
\begin{equation}\label{th1-eq1}
|y|_{C_\dbF([0,T];L^2(\Om;H))} +
|z|_{L^2_\dbF(0,T;\wt U)} \leq C(T)\big(
|y_0|_{L^2_{\cF_0}(\Om;H)} +
|u|_{L^2_\dbF(0,T;U)} \big).
\end{equation}
\end{theorem}
\begin{remark}
One can also consider stochastic Schr\"odinger
equations with variable coefficients. Following
the method in \cite{GS} and the proof of Theorem
\ref{th1}, one can see that Theorem \ref{th1} is
also true if the Laplacian operator in
\eqref{S-system1} is replaced by a general
elliptic operator. As we said before, to present
the key idea in a simple way, we do not pursue
the full technical generality.
\end{remark}
To prove Theorem \ref{th1}, we first write it in
an abstract form. To this end, let us define an
unbounded linear operator on $H$ as follows:
$$
\left\{
\begin{array}{ll}\ds
D(A)=H_0^1(G),\\
\ns\ds \langle Af,g
\rangle_{H^{-1}(G),H_0^1(G)}=\int_G \nabla
f(x)\cd\overline{\nabla g(x)}dx,\q \forall\,
f,g\in H_0^1(G).
\end{array}
\right.
$$

Define a map $\Upsilon:L^2(\G_0)\to L^2(G)$ as
follows:
$$
\Upsilon u =v,
$$
where $v$ is the solution to
$$
\left\{
\begin{array}{ll}\ds
\D v = 0 &\mbox{ in } G,\\
\ns\ds v=u &\mbox{ on } \G.
\end{array}
\right.
$$
From the definition of $\Upsilon$ and the
classical theory for elliptic equations with
non-homogeneous boundary condition(see
\cite[Chapter 2]{LM} for example), we get that
there is a constant $C>0$ such that for any
$(t,\om)\in (0,T)\times\Om$,
$$
\int_G |i\Upsilon u(t,\om)|^2dx\leq
C\int_{\G}|u(t,\om)|^2d\G,
$$
which deduces that
\begin{equation}\label{9.18-eq3}
\int_0^T\int_G |i\Upsilon u|^2dxdt \leq
C\int_0^T\int_{\G}|u|^2d\G dt.
\end{equation}

Define two operators $J,K\in
\cL(L^2_\dbF(0,T;H))$ as
$$
J h = a h,\q K h=bh,\q\forall\, h\in
L^2_\dbF(0,T;H).
$$
Then, the system \eqref{S-system1} can be
written as
\begin{equation}\label{S-system2}
dy-iA(y-\Upsilon u)=Jydt + KydW \qq\mbox{ in
}(0,+\infty).
\end{equation}

Clearly,
$$
D(A)\subset D(A^{\frac{1}{2}}) \hookrightarrow H
\hookrightarrow [D(A)^{\frac{1}{2}}]' \subset
[D(A)]'.
$$
Denote by $\wt A$ the extension of $A$ as a
bounded linear operator from
$D(A^{\frac{1}{2}})$ to $[D(A)^{\frac{1}{2}}]'$
as follows:
$$
\langle \wt Af,g
\rangle_{[D(A)^{\frac{1}{2}}]',D(A^{\frac{1}{2}})}
= \langle A^{\frac{1}{2}}f,
A^{\frac{1}{2}}g\rangle_{H},\q\forall\, f,g\in
D(A)^{\frac{1}{2}}.
$$
Then $i\wt A$ generates a $C_0$-group on
$[D(A)^{\frac{1}{2}}]'$. Thus, \eqref{S-system2}
can be reduced to
\begin{equation}\label{S-system3}
dy=i\wt A y+ Budt + Jydt + KydW \qq\mbox{ in
}(0,\infty),
\end{equation}
where $B\in \cL(U,[D(A)^{\frac{1}{2}}]')$ such
that
\begin{equation}\label{9.18-eq1}
B\eta=-i\wt A\Upsilon \eta,\q\forall\, \eta\in
U.
\end{equation}
Denote by $B^*$ the adjoint operator of $B$.
Then, $B^*\in \cL(D(A)^{\frac{1}{2}},U)$ and
$$
\langle B^* f,\eta \rangle_U = \langle f,B\eta
\rangle_{D(A)^{\frac{1}{2}},[D(A)^{\frac{1}{2}}]'},\q
\forall\,f\in D(A^{\frac{1}{2}}),\, \eta\in U.
$$
For any $f\in D(A)$ and $\eta\in
C_0^\infty(\G)$, we have that
$$
\begin{array}{ll}\ds
\langle f,B\eta
\rangle_{D(A)^{\frac{1}{2}},[D(A)^{\frac{1}{2}}]'}
\3n&\ds= \langle Af, \wt A^{-1}B\eta \rangle_H =
i\langle Af,  \Upsilon \eta \rangle_H = i\langle
A^{\frac{1}{2}}f, A^{-\frac{1}{2}}\Upsilon \eta
\rangle_{L^2(G)}\\
\ns&\ds = i\langle AA^{-1}f, \Upsilon \eta
\rangle_{L^2(G)}= -\Big\langle
i\frac{\pa((-\D)^{-1}f)}{\pa\nu}, \eta
\Big\rangle_{U}.
\end{array}
$$
Since $C_0^\infty(\G_0)$ is dense in
$L^2(\G_0)$, we have that
\begin{equation}\label{9.18-eq2}
B^* =
-i\frac{\pa((-\D)^{-1}f)}{\pa\nu}\Big|_{\G_0}.
\end{equation}
The system \eqref{S-system1} can be written in
the following abstract form:
\begin{equation}\label{S-system4}
\left\{
\begin{array}{ll}\ds
dy = i\wt Aydt + Bu dt + J y dt + K y dW(t)
&\mbox{ in } (0,+\infty),\\
\ns\ds y(0)=y_0,\\
\ns\ds \f=B^* y &\mbox{ in } (0,+\infty).
\end{array}
\right.
\end{equation}

Next, let us give the following result.
\begin{proposition}\label{9.18-prop1}
Let $\mu = \mu(x) =
(\mu^1,\cdots,\mu^d):\mathbb{R}^d \to
\mathbb{R}^d$ be a vector field of class $C^1$
and $\hat \f$ an
$H^2_{loc}(\mathbb{R}^d)$-valued $\dbF$-adapted
semi-martingale. Then for a.e. $x \in
\mathbb{R}^n$ and $\dbP$-a.s. $\omega \in
\Omega$, it holds that
\begin{equation}\label{9.18-prop1-eq1}
\begin{array}
{ll} & \ds \mu\cdot\nabla\bar{\hat \f}(d\hat \f
+ i\Delta \hat \f dt) -
\mu\cdot\nabla \hat \f(d\bar{\hat \f} - i\Delta \bar{\hat \f} dt)\\
\ns  =& \ds  \nabla\cd \Big[ i(\mu\cdot\nabla
\bar{\hat \f})\nabla \hat \f+ i(\mu\cdot\nabla
\hat \f)\nabla \bar{\hat \f}  - (\hat \f
d\bar{\hat \f}) \mu - i|\nabla
\hat \f|^2\mu \Big]dt + d(\mu\cd\nabla \bar{\hat \f} \hat \f)\\
\ns &  \ds - i\sum_{j,k=1}^d (\mu^k_j + \mu^j_k)
\hat \f_{j}\bar{\hat \f}_{k}dt + i(\nabla\cdot
\mu) |\nabla \hat \f|^2 dt + (\nabla\cdot \mu)
\hat \f d\bar{\hat \f} - (\mu\cd\nabla
d\bar{\hat \f}) d\hat\f.
\end{array}
\end{equation}
\end{proposition}
 {\it Proof of
Proposition \ref{9.18-prop1}} : The proof is a
direct computation. We have that
\begin{equation}\label{h1}
\begin{array}
{ll} & \ds i\sum_{k=1}^d\sum_{j=1}^d
\mu^k\bar{\hat \f}_k \hat \f_{jj}+
i\sum_{k=1}^d\sum_{j=1}^d \mu^k \hat \f_k \bar{\hat \f}_{jj}\\
\ns = & \ds
i\sum_{k=1}^d\sum_{j=1}^d\Big[(\mu^k\bar{\hat
\f}_k \hat \f_j)_j + (\mu^k \hat \f_k\bar{\hat
\f}_j)_j + \mu^k_k|\hat \f_j|^2 - (\mu^k|\hat
\f_j|^2)_k - (\mu^k_j +\mu^j_k ) \bar{\hat \f}_k
\hat \f_j \Big]
\end{array}
\end{equation}
and that
\begin{equation}\label{h2}
\begin{array}{ll}\ds
\q \sum_{k=1}^d(\mu^k\bar{\hat \f}_k d\hat \f-\mu^k \hat \f_k d\bar{\hat \f})\\
\ns\ds = \sum_{k=1}^d\Big[d(\mu^k\bar{\hat \f}_k
\hat \f) - \mu^k \hat \f d\bar{\hat \f}_k -
\mu^k d\bar{\hat \f}_k
d\hat \f -(\mu^k \hat \f d\bar{\hat \f})_k + \mu^k \hat \f d\bar{\hat \f}_k + \mu_k^k \hat \f d\bar{\hat \f} \Big]\\
\ns\ds = \sum_{k=1}^d\Big[d(\mu^k\bar{\hat \f}_k
\f) - \mu^k d\bar{\hat \f}_k d\hat \f -(\mu^k
\hat \f d\bar{\hat \f})_k + \mu_k^k \hat \f
d\bar{\hat \f} \Big].
\end{array}
\end{equation}
Combining \eqref{h1} and \eqref{h2}, we get the
equality \eqref{9.18-prop1-eq1}.
\endpf

\vspace{0.2cm}

Now we are in a position to prove Theorem
\ref{th1}.

{\it Proof of Theorem \ref{th1}}\,: We first
show that $B$ is an admissible control operator
with respect to the semigroup $\{S(t)\}_{t\geq
0}$ generated by $iA$, i.e., we prove that there
is a constant $C>0$ such that for any $u\in
L^2_\dbF(0,+\infty;U)$,
$$
|\Psi_T u|_{L^2_{\cF_T}(\Om;H)} \leq
C|u|_{L^2_\dbF(0,T;U)}
$$
To proof the above inequality, we only need to
establish the following inequality:
\begin{equation}\label{12.25-eq21}
|\Psi_T^* \xi|_{L^2_\dbF(0,T;U)} \leq
C|\xi|_{L^2_{\cF_T}(\Om;H)},\q \forall\,\xi\in
L^2_{\cF_T}(\Om;H),
\end{equation}
where $C=C(T)$ is independent of $\xi$. To
achieve this goal, we consider the following
backward stochastic Schr\"odinger equation:
\begin{equation}\label{ad-S-system1}
\left\{
\begin{array}{ll}\ds
dv + i\D v dt = -\bar av dt -\bar b V dt +
VdW(t)
&\mbox{ in } G\times (0,T),\\
\ns\ds v=0 &\mbox{ on }\G\times (0,T),\\
\ns\ds v(T)=v_T &\mbox{ in }G,
\end{array}
\right.
\end{equation}
where $v_T\in L^2_{\cF_T}(\Om;H^{-1}(G))$. By
the classical theory of backward stochastic
evolution equations(see \cite{HP}), we know that
the equation \eqref{ad-S-system1} admits a
unique solution $(v,V)\in
C_\dbF([0,T];L^2(\Om;H^{-1}(G)))\times
L^2_\dbF(0,T;H^{-1}(G))$.

Let $w=A^{-1}v$. Then $w$ solves that
\begin{equation}\label{S-system5}
\left\{
\begin{array}{ll}\ds
dw + i\D wdt=   - A^{-1}(\bar av)dt -
A^{-1}(\bar bV)dt + A^{-1}VdW(t) &\mbox{ in }
G\times
(0,T),\\
\ns\ds w=0 &\mbox{ on } \G\times (0,T),\\
\ns\ds w(T)=w_T= A^{-1}v_T &\mbox{ in } G,\\
\ns\ds B^* v = B^* AA^{-1}v = -i\frac{\pa
w}{\pa\nu} &\mbox{ on } \G_0\times (0,T).
\end{array}
\right.
\end{equation}
From the definition of $A$, we know that
$$
\begin{array}{ll}\ds
w\in C_\dbF([0,T];L^2(\Om;H_0^1(G))),\q
A^{-1}V\in
L^2_\dbF(0,T;H_0^1(G)),\\
\ns\ds A^{-1}(\bar av)\in
L^\infty_\dbF(0,T;H^1_0(G)),\q A^{-1}(\bar
bV)\in L^2_\dbF(0,T;H^1_0(G)).
\end{array}
$$
Further,
\begin{equation}\label{9.18-eq6}
\left\{
\begin{array}{ll}\ds
|w|_{C_\dbF([0,T];L^2(\Om;H_0^1(G)))}\leq
C|v|_{C_\dbF([0,T];L^2(\Om;H^{-1}(G)))},\\
\ns\ds |A^{-1}V|_{L^2_\dbF(0,T;H_0^1(G))}\leq
C|V|_{L^2_\dbF(0,T;H^{-1}(G))},\\
\ns\ds |A^{-1}(\bar
av)|_{L^\infty_\dbF(0,T;H^1_0(G))}\leq
C|a|_{L^\infty_\dbF(0,T;W_0^{1,\infty}(G))}|v|_{L^2_\dbF(0,T;H^{-1}(G))},\\
\ns\ds |A^{-1}(\bar
bV)|_{L^\infty_\dbF(0,T;H^1_0(G))}\leq
C|b|_{L^\infty_\dbF(0,T;W_0^{1,\infty}(G))}|V|_{L^2_\dbF(0,T;H^{-1}(G))}.
\end{array}
\right.
\end{equation}

Since $\G$ is $C^2$, there is a $C^1$ vector
field $h=(h^1,\cds,h^d):\overline G\to\dbR^d$
such that
$$
h(x)=\nu(x) \q\mbox{ on }\G,\q |h(x)|\leq
1,\q\forall\, x\in G.
$$
Let us take $\mu=h$ and $\hat \f = w$ in
\eqref{9.18-prop1-eq1}. Integrating it in
$G\times (0,T)$ and taking the mathematical
expectation, we have that
\begin{equation}\label{9.18-eq4}
\!\!\begin{array} {ll} \ds \q
-\mE\int_0^T\!\int_G\!
h\cdot\nabla\bar{w}\big[A^{-1}(\bar av)+
A^{-1}(\bar bV)\big]dxdt + \mE\int_0^T\!\!\int_G
h\cdot\nabla w\big[ A^{-1}(a\bar v) +
A^{-1}(b\overline V)\big]dxdt\\
\ns  = \ds  i\mE\int_0^T\int_{\G} \Big|\frac{\pa
w}{\pa\nu}\Big|^2dxdt
+ \int_G h\cd\nabla \bar{w}(T)w(T)dx - \int_G h\cd\nabla \bar{w}(0)w(0)dx\\
\ns \ds \q - i\mE\int_0^T\int_G\sum_{j,k=1}^d
(h^k_j + h^j_k) w_{j}\bar{w}_{k}dxdt +
i\mE\int_0^T\int_G(\nabla\cdot h) |\nabla w|^2
dxdt\\
\ns\ds \q + \mE\int_0^T\!\int_G(\nabla\cdot h) w
\big[ i\D\bar w -\! A^{-1}(a\bar v)dt -\!
A^{-1}(b\overline V)dt \big]dxdt -
\mE\int_0^T\!\int_G\big[h\cd\nabla
A^{-1}(\overline{V})\big] A^{-1}Zdxdt.
\end{array}
\end{equation}
From \eqref{9.18-eq6}, we know that
\begin{equation}\label{9.18-eq5}
\begin{array}{ll}\ds
\q\Big|\mE\int_0^T\int_G
h\cdot\nabla\bar{w}\big[A^{-1}(\bar av) +
A^{-1}(\bar bV)\big]dxdt\Big| +
\Big|\mE\int_0^T\int_G h\cdot\nabla w\big[
A^{-1}(a\bar v) +
A^{-1}(b\overline V)\big]dxdt\Big|\\
\ns\ds \leq C\big(|w|^2_{L^2_\dbF(0,T;H_0^1(G))}
+
|a|^2_{L^\infty_\dbF(0,T;W_0^{1,\infty}(G))}|v|^2_{L^2_\dbF(0,T;H^{-1}(G))}
+
|b|^2_{L^\infty_\dbF(0,T;W_0^{1,\infty}(G))}|V|^2_{L^2_\dbF(0,T;H^{-1}(G))}
\big)\\
\ns\ds \leq
C\big(|v|^2_{C_\dbF([0,T];L^2(\Om;H^{-1}(G)))} +
|a|^2_{L^\infty_\dbF(0,T;W_0^{1,\infty}(G))}|v|^2_{L^2_\dbF(0,T;H^{-1}(G))}
\\
\ns\ds \qq +
|b|^2_{L^\infty_\dbF(0,T;W_0^{1,\infty}(G))}|V|^2_{L^2_\dbF(0,T;H^{-1}(G))}
\big)\\
\ns\ds \leq C
|v_T|^2_{L^2_{\cF_T}(\Om;H^{-1}(G))},
\end{array}
\end{equation}
\begin{equation}\label{9.18-eq7}
\begin{array}{ll}\ds
\q\Big|\int_G h\cd\nabla \bar{w}(T)w(T)dx -
\int_G h\cd\nabla \bar{w}(0)w(0)dx\Big|\\
\ns\ds \leq
C|w|^2_{C_\dbF([0,T];L^2(\Om;H_0^1(G)))}\leq
C|v|^2_{C_\dbF([0,T];L^2(\Om;H^{-1}(G)))}\leq C
|v_T|^2_{L^2_{\cF_T}(\Om;H^{-1}(G))},
\end{array}
\end{equation}
\begin{equation}\label{9.18-eq8}
\begin{array}{ll}\ds
\q\Big|- i\mE\int_0^T\int_G\sum_{j,k=1}^d
\big(h^k_j + h^j_k\big) w_{j}\bar{w}_{k}dxdt +
i\mE\int_0^T\int_G(\nabla\cdot h) |\nabla w|^2
dxdt\Big|\\
\ns\ds \leq C|w|^2_{L^2_\dbF(0,T;H_0^1(G))}\leq
C|w|^2_{C_\dbF([0,T];L^2(\Om;H_0^1(G)))}\leq
C|z|^2_{C_\dbF([0,T];L^2(\Om;H^{-1}(G)))}\leq C
|z_T|^2_{L^2_{\cF_T}(\Om;H^{-1}(G))},
\end{array}
\end{equation}
and
\begin{equation}\label{9.18-eq9}
\begin{array}{ll}\ds
\q\Big|\mE\int_0^T\int_G\big[h\cd\nabla
A^{-1}(\overline{V})\big] A^{-1}Vdxdt\Big|\\
\ns\ds \leq
C|A^{-1}V|_{L^2_\dbF(0,T;H_0^1(G))}|V|_{L^2_\dbF(0,T;L^2(G))}\leq
C|A^{-1}V|^2_{L^2_\dbF(0,T;H^{-1}(G))}\leq C
|v_T|^2_{L^2_{\cF_T}(\Om;H^{-1}(G))}.
\end{array}
\end{equation}
Further, since
$$
\begin{array}{ll}\ds
\q\mE\int_0^T\int_G(\nabla\cdot h) w  i\D\bar
wdxdt\\
\ns \ds =
-\mE\int_0^T\int_G\nabla\cd(\nabla\cdot h) w
i\nabla\bar wdxdt -\mE\int_0^T\int_G(\nabla\cdot
h) \nabla w i\nabla\bar wdxdt,
\end{array}
$$
due to \eqref{9.18-eq6}, we have that
\begin{equation}\label{9.18-eq10}
\begin{array}{ll}\ds
\q\Big|\mE\int_0^T\int_G(\nabla\cdot h) w \big[
i\D\bar w - A^{-1}(a\bar v)dt -
A^{-1}(b\overline V)dt
\big]dxdt\Big|\\
\ns\ds \leq C\big(|w|^2_{L^2_\dbF(0,T;H_0^1(G))}
+
|a|^2_{L^\infty_\dbF(0,T;W_0^{1,\infty}(G))}|v|^2_{L^2_\dbF(0,T;H^{-1}(G))}
+
|b|^2_{L^\infty_\dbF(0,T;W_0^{1,\infty}(G))}|V|^2_{L^2_\dbF(0,T;H^{-1}(G))}\big)\\
\ns\ds \leq C
|v_T|^2_{L^2_{\cF_T}(\Om;H^{-1}(G))}.
\end{array}
\end{equation}
From \eqref{9.18-eq4} to \eqref{9.18-eq10}, we
obtain that
$$
\mE\int_0^T\int_{\G} \Big|\frac{\pa
w}{\pa\nu}\Big|^2dxdt \leq
C|z_T|^2_{L^2_{\cF_T}(\Om;H^{-1}(G))},
$$
which implies that
\begin{equation}\label{9.18-eq11}
\mE\int_0^T\int_{\G} \big|B^* z\big|^2dxdt \leq
C|z_T|^2_{L^2_{\cF_T}(\Om;H^{-1}(G))}.
\end{equation}
Hence, we know that  $B$ is an admissible
control operator.

Next, we prove that the input/output map is a
bounded linear operator. Let $\tilde w=A^{-1}y$.
Then $\tilde w$ solves
\begin{equation}\label{S-system6}
\left\{
\begin{array}{ll}\ds
d\tilde w + i\D \tilde wdt= -i\Upsilon u dt +
A^{-1}(ay)dt - A^{-1}(by)dW(t)  &\mbox{ in }
G\times
(0,+\infty),\\
\ns\ds \tilde w=0 &\mbox{ on } \G\times (0,+\infty),\\
\ns\ds \tilde w(0)=\tilde w_0= A^{-1}y_0 &\mbox{ in } G,\\
\ns\ds \f=-i\frac{\pa \tilde w}{\pa\nu} &\mbox{
on } \G_0\times (0,+\infty).
\end{array}
\right.
\end{equation}
Since $B$ is a admissible control operator, we
know that $y\in C_\dbF([0,T];L^2(\Om;H))$ and
\begin{equation}\label{9.18-eq12}
|y|_{C_\dbF([0,T];L^2(\Om;H))}\leq
C\big(|y_0|_{L^2_{\cF_0}(\Om;H)}+|u|_{L^2_\dbF(0,T;U)}\big).
\end{equation}

From the definition of $A$, we know that
$$
\begin{array}{ll}\ds
\tilde w\in C_\dbF([0,T];L^2(\Om;H_0^1(G))), \q
A^{-1}(ay)\in L^\infty_\dbF(0,T;H^1_0(G)),\q
A^{-1}(by)\in L^\infty_\dbF(0,T;H^1_0(G)).
\end{array}
$$
Further, thanks to \eqref{9.18-eq12}, we have
that
$$
|\tilde w|_{C_\dbF([0,T];L^2(\Om;H_0^1(G)))}\leq
C|y|_{C_\dbF([0,T];L^2(\Om;H^{-1}(G)))}\leq
C\big(|y_0|_{L^2_{\cF_0}(\Om;H)}+|u|_{L^2_\dbF(0,T;U)}\big),
$$
$$
\begin{array}{ll}\ds
|A^{-1}(ay)|_{L^\infty_\dbF(0,T;H^1_0(G))}\3n&\ds\leq
C|a|_{L^\infty_\dbF(0,T;W_0^{1,\infty}(G))}|y|_{L^2_\dbF(0,T;H^{-1}(G))}\\
\ns&\ds\leq
C\big(|y_0|_{L^2_{\cF_0}(\Om;H)}+|u|_{L^2_\dbF(0,T;U)}\big)
\end{array}
$$
and
$$
\begin{array}{ll}\ds
|A^{-1}(by)|_{L^\infty_\dbF(0,T;H^1_0(G))}\3n&\ds\leq
C|b|_{L^\infty_\dbF(0,T;W_0^{1,\infty}(G))}|y|_{L^2_\dbF(0,T;H^{-1}(G))}\\
\ns&\ds\leq
C\big(|y_0|_{L^2_{\cF_0}(\Om;H)}+|u|_{L^2_\dbF(0,T;U)}\big).
\end{array}
$$
According to these inequalities, similar to the
proof of \eqref{9.18-eq11}, we can obtain that
\begin{equation}\label{9.18-eq14}
\mE\int_0^T\int_{\G} \Big|\frac{\pa
(-\D)^{-1}y}{\pa\nu}\Big|^2dxdt \leq
\mE\int_0^T\int_{\G} \Big|\frac{\pa \tilde
w}{\pa\nu}\Big|^2dxdt \leq
C\big(|y_0|^2_{L^2_{\cF_0}(\Om;H)}+|u|^2_{L^2_\dbF(0,T;U)}\big).
\end{equation}
This implies that the boundary observation
operator in the system \eqref{S-system1} is
admissible and the system \eqref{S-system1} is
well-posed.
\endpf


\section{Further comments and open
problems}\label{sec-5}


This paper is only a first and basic attempt to
the study of well-posed linear stochastic
system. In my opinion, there are many
interesting and important problems in this
topic. We present some of them here briefly:

\begin{itemize}

\item{\bf Further properties for stochastic well-posed linear systems.}

Deterministic well-posed linear systems enjoy
many deep and useful properties(see
\cite{Staffans} for example). In this paper, we
only investigate some very basic ones. It is an
interesting and maybe challenging problem  to
study what kind of properties for deterministic
well-posed linear systems also holds for
stochastic well-posed linear systems.

\medskip

\item {\bf The well-posedness of stochastic partial
differential equations with boundary
control/observation}.

The main motivation of introducing stochastic
well-posed linear systems is to study the
stochastic partial differential equations with
boundary control/observation. In this paper, we
only consider two special examples. It is well
known that many deterministic partial
differential equations with boundary
control/observation are well-posed(see
\cite{GS,GZ1,GZ2} and the rich references
therein). It deserves to generalize there
results for stochastic partial differential
equations.

\medskip

\item {\bf The study of the stochastic regular system}.

In the deterministic framework, there is another
important concept, i.e., regular systems, in the
study of infinite dimensional linear systems,
which has very close relation to the well-posed
system. One can also define the stochastic
regular systems and study their properties. Some
of them will appear in our forthcoming paper
\cite{Lu1}. But there are still lots of problems
should be studied.

\medskip

\item {\bf The stabilization of stochastic control systems}.

Once we prove that a system is well-posed with
$U=\wt U$, then we can consider the
stabilization of the system by the state
feedback control. Such kind of problems are
extensively studied for the deterministic
control systems in the literature(see
\cite{Logemann,LRT,Mikkola,WR} and the rich
references therein). On the other hand, as far
as we know, there is no result for the
stochastic counterpart.

\end{itemize}

\end{document}